\documentclass{amsart}
\usepackage{amsmath}
\usepackage{amsfonts}
\usepackage{amssymb}
\usepackage{graphicx}
\usepackage{multirow}

\numberwithin{equation}{section}
\newcommand\T{\rule{0pt}{2.5ex}}
\begin{document}

\title{On centers and central lines of triangles \mbox {in the elliptic plane}}
\author{Manfred Evers}
\curraddr[Manfred Evers]{Bendenkamp 21, 40880 Ratingen, Germany}
\email[Manfred Evers]{manfred\_evers@yahoo.com}
\date{\today}

\begin{abstract}
We determine barycentric coordinates of triangle centers in the elliptic plane. The main focus is put on centers that lie on lines whose euclidean limit (triangle excess $\rightarrow  0$) is the Euler line or the Brocard line. We also investigate curves which can serve in elliptic geometry as substitutes for the euclidean nine-point-circle, the first Lemoine circle or the apollonian circles.
\end{abstract}

\maketitle

\section*{Introduction}\hspace*{\fill} \\
\noindent In the first section we give a short introduction to metric geometry in the projective plane. We assume the reader is familiar with this subject, but we recall some fundamental definitions and theorems, in order to introduce the terminology and to fix notations. 
The second section provides appropriate tools (definitions, theorems, rules) for calculating the barycentric coordinates (in section 3) of a series of centers lying on four central lines of a triangle in the elliptic plane.\\ 
The content of this work is linked to results presented by  Wildberger \cite{27}, Wildberger and Alkhaldi [28], Ungar [25], Horv\'ath [8], Vigara [26], Russell [18]. \vspace*{1.5 mm}

\section{Metric geometry in the projective plane \vspace*{0.5 mm}}

\subsection{} \textbf{The projective plane, its points and its lines}\hspace*{\fill} \\
\noindent\hspace*{5mm}Let $V$ be the three dimensional vector space $\mathbb{R}^{3}$, equipped with the canonical dot product $\boldsymbol{p} \cdot \boldsymbol{q} = (p_0,p_1,p_2) \cdot (q_0,q_1,q_2) = p_0q_0+p_1q_1+p_2q_2$ and the induced norm $\|.\|$, and let
$\mathcal{P}$ denote the projective plane $({V}{-}\{\boldsymbol{0}\})/\mathbb{R}^\times$. The image of a non-zero vector $\boldsymbol{p} = (p_0,p_1,p_2) \in V $ under the canonical projection $\Pi : V \rightarrow \mathcal{P}$ will be denoted by $(p_0:p_1:p_2)$ and will be regarded as a point in this plane. \\
\noindent\hspace*{5mm}Given two different points  $P$ and $Q$ in this projective plane, there exists exactly one line that is incident with these two points. It is called the \textit{join} $P \vee Q$ of $P$ and $Q$.  If $\boldsymbol{p} = (p_0,p_1,p_3)$, $\boldsymbol{q} = (q_0,q_1,q_2)$ are two non zero vectors with $\Pi(\boldsymbol{p}) = P$ and $\Pi(\boldsymbol{q}) = Q$, then the line $P \vee Q$ through $P$ and $Q$ is the set of points $\Pi(s \boldsymbol{p} + t \boldsymbol{q})$ with $s, t \in \mathbb{R}$.
One can find linear forms $l \in V^{*}{-}\{\boldsymbol{0^*}\}$ with $ker(l) = span(\boldsymbol{p}, \boldsymbol{q})$. A suitable $l$ is, for example, $l = \ast(\boldsymbol{p}\times \boldsymbol{q}) = (\boldsymbol{p}\times \boldsymbol{q})^*$, where $\times$ stands for the canonical cross product on $V = \mathbb{R}^{3}$ and $\ast$ for the isomorphism $V \rightarrow V^*, \ast(\boldsymbol{r}) = (.){\cdot}\boldsymbol{r}.$ The linear form $l$ is uniquely determined up to a nonzero real factor, so there is a ${1{:}1}$-correspondence between the lines in the projective plane and the elements of $\mathcal{P}^* = (V^{*}{-}\{\boldsymbol{0^*}\})/\mathbb{R}^\times$. We identify the line $l = P \vee Q$ with the element $(p_1q_2 - p_2q_1:p_2q_0 - p_0q_2:p_0q_1 - p_1q_0)^* \in \mathcal{P}^*.$\\
\noindent\hspace*{5mm}In the projective plane, two different lines $k = (k_0:k_1:k_2)^*, l = (l_0:l_1:l_2)^*$ always meet in one point $k \wedge l = \Pi((k_0,k_1,k_2)\times (l_0,l_1,l_2))$, the so called \textit{meet} of these lines.

\subsection{} \textbf{Visualizing points and lines}\hspace*{\fill} \\
\noindent\hspace*{5mm}%
Using an orthogonal coordinate system, we know how to visualize (in a canonical way) a point with cartesian coordinates $(p_1,p_2)$ in the affine plane $\mathbb{A}^2$.  A point $P = (1:p_1:p_2)\in \mathcal{P}$ will be visualized as the point $(p_1,p_2)$ in the affine plane, and we will call $P^\vee := (1,p_1,p_2) \in \mathbb{R}^{3}$ the visualizing vector of $P$.  But in $\mathcal{P}$ there exist points $(0:p_1:p_2)$ which can not be visualized in the affine plane. These points are considered to be points on the ``line at infinity''. For these points we define $P^\vee := (0,1,p_2/p_1)$, if $p_1 \ne 0$, and $P^\vee := (0,0,1)$, otherwise. In this way, we ensure that the triple $P^\vee$ is strictly positive with respect to the lexicographic order.

A line appears as a ``stright line'' in the coordinate system.

\subsection{} \textbf{Collineations and correlations}\hspace*{\fill} \\
\noindent\hspace*{5mm}%
A \textit{collineation} on $\mathcal{P}$ is a bijective mapping $\mathcal{P} \rightarrow \mathcal{P}$ that maps lines to lines. These collineations form the group of automorphisms of $\mathcal{P}$. \\Collineations preserve the cross ratio $(P,Q;R,S)$ of four points on a line.

A \textit{correlation} on the projective plane is either a point-to-line transformation that maps collinear points to concurrent lines, or it is a line-to-point transformation that maps concurrent lines to collinear points.\\  Correlations, as collineationes, preserve the cross ratio.

\subsection{} \textbf{Metrical structures on  $\mathcal{P}$}
\subsubsection{} \textbf{The absolute conic }\hspace*{\fill} \\
\noindent\hspace*{5mm}%
One of the correlations is the \textit{polarity} with respect to the \textit{absolute conic} $\mathcal{C}$. This correlation assigns each point $P = (p_0:p_1:p_2)$ its polar line $P^\delta = (p_0:\sigma p_1: \varepsilon \sigma p_2)^*$ and assigns each line $l = (l_0:l_1:l_2)^*$ the corresponding pole $l^\delta  = (\varepsilon \sigma l_0: \varepsilon l_1:l_2)$. Here, $\sigma \in \mathbb{R^\times}$, $\varepsilon \in \{-1,1\}$,  and the absolute conic consists of all points $P$ with $P \in P^\delta$. $P^\delta$ and $l^\delta$ are called the dual of $P$ and $l$, respectively.

Besides the norm $\|.\|$ we introduce a seminorm $\|.\|_{\sigma,\varepsilon}$ on $V$; this is defined by: $\|\boldsymbol{p}\|_{\sigma,\varepsilon} = \sqrt{|p_0^2 + \sigma p_1^2 + \varepsilon \sigma  p_2^2|}$ for $\boldsymbol{p} = (p_0,p_1,p_2)$. It can be easily checked that $P = \Pi(\boldsymbol{p})\in \mathcal{C}$ exactly when $\|\boldsymbol{p}\|_{\sigma,\varepsilon} = 0$. Points on $\mathcal{C}$ are called \textit{isotropic points}. 

The dual of an isotropic point is an \textit{isotropic line}. As isotropic points, also isotropic lines form a conic, the dual conic $\mathcal{C^\delta}$ of $\mathcal{C}$. 

\subsubsection{} \textbf{Cayley-Klein geometries}\hspace*{\fill} \\
\noindent\hspace*{5mm}%
Laguerre and Cayley were presumably the first to recognize that conic sections can be used to define the angle between lines and the distance between points, cf. \cite{1}. An important role within this connection plays the cross ratio of points and of lines. We do not go into the relationship between conics and measures; there are many books and articles about Cayley-Klein-geometries (for example [13, 15, 17, 22, 23]) treating this subject. Particularly extensive investigations on cross-ratios offers Vigara [26].

Later, systematic studies by Felix Klein [13] led to a classification of metric geometries on $\mathcal{P}.$ He realized that not only a geometry determines its automorphisms, but one can make use of automorphisms to define a geometry [12]. The automorphisms on $\mathcal{P}$ are the collineations. By studying the subgroup of collineations that keep the absolute conic fixed (as a whole, not pointwise), he was able to find different metric geometries on $\mathcal{P}.$  

If $\varepsilon = 1$ and $\sigma > 0$, there are no real points on the absolute conic, so there are no isotropic points and no isotropic lines. The resulting geometry was called \textit{elliptic} by Klein. It is closely related to spherical geometry. From the geometry on a sphere we get an elliptic geometry by identifying antipodal points. Already Riemann had used spherical geometry to get a new metric geometry with constant positive curvature, cf. [14, ch. 38].

Klein [13] also showed that in the elliptic case the euclidean geometry can be received as a limit for $\sigma \rightarrow 0$. (For $\sigma \rightarrow \infty$ one gets the polar-euclidean geometry.)

\subsection{} \textbf{Metrical structures in the elliptic plane}\hspace*{\fill}\\
\noindent\hspace*{5mm}%
In the following, we consider just the elliptic case. Thus, we assume $\varepsilon$ = $1$ and $\sigma > 0$. Nearly all our results can be transferred to other Cayley-Klein geometries, some even to a ``mixed case'' where the points lie in different connected components of $\mathcal{P} -\mathcal{C}$, cf. [10, 18, 25, 27, 28]. Nevertheless, it is less complicated to derive results in the elliptic case, because: First, there are no isotropic points and lines. Secondly, if we additionally put $\sigma$ to $1$ - and this is what we are going to do - , then the norm $\|.\|_{\sigma,\varepsilon}$ agrees with the standard norm $\|.\|$ and this simplifies many formulas. 
For example, we have $(p_0{:}p_1{:}p_2)^\delta = (p_0{:}p_1{:}p_2)^*.$

\subsubsection{} \textbf{Barycentric coordinates of points}\hspace*{\fill} \\
\noindent\hspace*{5mm}%
For $P \in \mathcal{P}$, define the vector $P^\circ$ by $P^\circ:= P^\vee/\|P^\vee\|$. Given n points $P_1,\cdots,P_n \in \mathcal{P}$, we say that a point $P$ is a (\textit{linear}) \textit{combination} of $P_1,\cdots,P_n$, if there are real numbers $t_1,\cdots,t_n$ such that $P = \Pi(t_1 P_1^\circ + \cdots + t_n P_n^\circ)$, and we write $P = t_1 P_1 + \cdots + t_n P_n$. \\
The points $P_1,\cdots,P_n$ form a \textit{dependent system} if one of the $n$ points is a combination of the others. Otherwise, $P_1,\cdots,P_n$ are \textit{independent}. A single point is always independent, so are two different points. Three points are independent exactly when they are not collinear. And more than three points in $\mathcal{P}$ always form a dependent system.

If $\Delta = ABC$ is a triple of three non-collinear points $A, B, C$, then every point $P \in \mathcal{P}$ can be written as a combination of these. If $P = s_1 A + s_2 B + s_3 C$ and $P = t_1 A + t_2 B + t_3 C$ are two such combinations, then there is always a real number $c \ne 0$ such that $c(s_1,s_2,s_3) = (t_1,t_2,t_3)$. Thus, the point $P$ is determined by $\Delta$ and the homogenous triple $(s_1:s_2:s_2)$. We write $P = [s_1:s_2:s_2]_\Delta$ and call this the \textit{representation of} $P$ \textit{by barycentric coordinates} with respect to $\Delta$. The terminology is not uniform here; the coordinates are also named \textit{gyrobarycentric} (Ungar [25]), \textit{circumlinear} (Wildberger, Alkhaldi [28]), \textit{triangular} (Horv\'ath [10]). 

Barycentric coordinates of a point $P$ can be calculated as follows: Because $A^\circ,B^\circ,C^\circ$ form a basis of $\mathbb{R}^{3}$, there is a unique way of representing $P^\circ$ by a linear combination $P^\circ = s_1 A^\circ + s_2 B^\circ + s_3 C^\circ$ of the base vectors; and the coordinates are 
\begin{equation*}
 {{s_1 = \frac{P^\circ \cdot (B^\circ \times C^\circ)}{A^\circ \cdot (B^\circ \times C^\circ)}},{\,s_2 = \frac{P^\circ \cdot (C^\circ \times A^\circ)}{B^\circ \cdot (C^\circ \times A^\circ)}},{\,s_3 = \frac{P^\circ \cdot (A^\circ \times B^\circ)}{C^\circ \cdot (A^\circ \times B^\circ)}}.} 
\end{equation*}\\
Since $A^\circ\!\cdot\!(B^\circ \times C^\circ) = B^\circ\!\cdot\!(C^\circ \times A^\circ) =$ ${C^\circ\!\cdot\!(A^\circ \times B^\circ)}$, we get 
\begin{equation*}
{s_1:s_2:s_3 = P^\circ\!\cdot\!(B^\circ \times C^\circ) : P^\circ\!\cdot\!(C^\circ \times A^\circ) : P^\circ\!\cdot\!(A^\circ \times B^\circ).}
\end{equation*}

\subsubsection{} \textbf{Orthogonality}\hspace*{\fill} \\
\noindent\hspace*{5mm}
We define orthogonality via polarity: A line $k$ is \textit{orthogonal} (or \textit{perpendicular}) to a line $l$ exactly when the dual $k^\delta$ of $k$ is a point on  $l$. It can easily be shown that, if $k$ is orthogonal to $l$, then $l$ is orthogonal to $k$.\\ 
The orthogonality between points is also defined: Two points are orthogonal precisely when their dual lines are. We can make use of the dot product to check if two points are orthogonal: Two points $P$ and $Q$ are orthogonal precisely when $P^\circ{\cdot}\,Q^\circ = 0$. 
Obviously, the set of points orthogonal to a point $P$ is its polar line $P^\delta$. 

\subsubsection{} \textbf{The distance between points and the length of line segments\\\hspace*{10 mm} in elliptic geometry}\hspace*{\fill} \\
\noindent\hspace*{5mm}%
Lines in elliptic geometry are without boundary. They are all the same length, usually set to $\pi$; this equals one half of the length of the great circle on a unit sphere $\mathbb{S}^2 \subset \mathbb{R}^3$ . In this case, the distance between two points $P$ and $Q$ is $d(P,Q) = \varphi$ with $\cos(\varphi) = |P^\circ{\cdot}\,Q^\circ|$. $P^\circ{\cdot}\,Q^\circ$ always takes values in the interval $]-1;1]$, so the distance between two points $P$ and $Q$ lies in the interval $[0,\frac {\pi}{2}]$, and $d(P,Q) = \frac {\pi}{2}$ implies that $P^\circ{\cdot}\,Q^\circ = 0$ and $P$ and $Q$ are orthogonal.

Two different points $P$ and $Q$ determine the line $P \vee Q$. The set $P \vee Q - \{P, Q\}$ consists of two connected components, the closure of these are called the \textit{line segments} of ${P,Q}$. One of these two segments contains all points $\Pi(s P^\circ+(1-s) Q^\circ)$ with $s(1-s){\ge}0$, while the other contains all points $\Pi(s P^\circ+(1-s) Q^\circ)$ with $s(1-s){\le}0$. The first segment will be denoted by $[P,Q]_{+}$, the second by $[P,Q]_{-}$. 

We show that $P+Q$ is the midpoint of $[P,Q]_+$ by proving the equation $(P+Q)^\circ{\cdot}P^\circ = (P+Q)^\circ{\cdot}Q^\circ$:

\centerline{$(P+Q)^\circ{\cdot}P^\circ = \frac{(P^\circ+Q^\circ)\cdot P^\circ}{\sqrt{(P^\circ+Q^\circ)\cdot(P^\circ+Q^\circ)}} = \frac {1+P^\circ\cdot Q^\circ}{\sqrt{2(1+P^\circ\cdot Q^\circ)}} = \sqrt{\frac{1+P^\circ\cdot Q^\circ}{2}} = (P+Q)^\circ{\cdot}Q^\circ$.} 
\noindent In the same way it can be verified that $P-Q$ is the midpoint of $[P,Q]_-$. Since $(P^\circ+Q^\circ){\cdot}(P^\circ - Q^\circ) = 0$, the two points $P+Q$ and $P-Q$ are orthogonal.\\
We now can calculate the measures (lengths) of the segments $[P,Q]_+$ and $[P,Q]_-$: \\
\centerline{$\mu([P,Q]_+) = \arccos(P^\circ \cdot Q^\circ) \in [0,\pi[$ and $\mu([P,Q]_-) = \pi - \mu([P,Q]_+)$.} \\
For further calculations the following formula will be useful:\\
\centerline{$ \sin(\mu([P,Q]_+)) = \sin(\mu([P,Q]_-)) = \sin(d(P,Q)) = \|P^\circ \times Q^\circ\|.$}
\noindent {\textit{Proof} of this formula}: $\sin(\mu([P,Q]_+)) = \sin(\mu([P,Q]_-)) = \sin(d(P,Q))$, because $\sin(\pi-x) = \sin(x)$ for $x\in [0,\pi[.$
The correctness of the last equation can be proved by verifying the equation $\|P^\circ{\times}\,Q^\circ\|^2 = 1 - (P^\circ{\cdot}\,Q^\circ)^2.\;\;\Box$

\subsubsection{} \textbf{Angles}\hspace*{\fill} \\
\noindent\hspace*{5mm}%
The (angle) distance between two lines $k$ and $l$ we get by dualizing the distance of two points: $d(k,l) = d(k^\delta,l^\delta)$. We even use the same symbol $d$ for the distance between lines as between points and do not introduce a new sign.\\  
By dualizing line segments, we get angles as subsets of the pencil of lines through a point which is the vertex of this angle: Given three different points $Q, R$ and $S$, we define the angles 
\begin{equation*}
\angle_{+} QSR := \{S \vee P |\, P \in [Q,R]_+\} \;\; \text{and} \;\;\angle_{-} QSR := \{S \vee P |\, P \in [Q,R]_-\}.
\end{equation*}
Using the same symbol $\mu$ for the measure of angles as for line segments, we have
\begin{equation*}
\mu(\angle_{+} QSR) = \arccos \frac{(S^\circ{\times}Q^\circ)}{\|(S^\circ{\times}Q^\circ)\|} \cdot \frac{(S^\circ{\times}R^\circ)}{\|(S^\circ{\times}R^\circ)\|} \;\;\text{and}\;\; \mu(\angle_{-} QSR) = \pi - \mu(\angle_{+} QSR).
\end{equation*}

\subsubsection{} \textbf{Perpendicular line through a point/perpendicular point on a line}\hspace*{\fill} \\
\noindent\hspace*{5mm}%
Consider a point $P = (p_0:p_1:p_2)$ and a line $l = (l_0:l_1:l_2)^*.$ We assume that $P \ne l^\delta.$  
The perpendicular from $P$ to $l$ is \begin{equation*}
\textrm{perp}(l,P) := P \vee l^\delta = (p_1l_2{-}p_2l_1: p_2l_0{-}p_0l_2: p_0l_1{-}p_1l_0)^*. 
\end{equation*}                                              
The line $\textrm{perp}(l,P)$ intersects $l$ at the point
\begin{equation*}
\begin{split}
Q &= l\wedge \textrm{perp}(l,P) \\
=&(l_0(l_1 p_1{+}l_2 p_2){-}p_0 (l_1^2{+}l_2^2): l_1 (l_0 p_0{+}l_2 p_2){-}p_1 (l_0^2{+}l_2^2): l_2 (l_0 p_0{+}l_1 p_1){-}p_2 (l_0^2{+}l_1^2)).
\end{split}
\end{equation*}
This point $Q$ is called the \textit{orthogonal projection} of $P$ on $l$ or the \textit{pedal} of $P$ on $l$.\hspace*{\fill} \\

Given two different points $P$ and $Q$, the perpendicular bisector of $[P,Q]_+$ is the line $(P+Q)\vee (P \vee Q)^\delta$ and the perpendicular bisector of $[P,Q]_-$ is the line $(P-Q)\vee (P \vee Q)^\delta$. A point on either of these perpendicular bisectors has the same distance from the endpoints
$P$ and $Q$ of the segment. \\

There is exactly one point $Q$ on the line $l$ with $d(Q,P) = \pi/2$ ; this is 
\begin{equation*}
Q = l \wedge P^\delta = (p_1l_2{-}p_2l_1: p_2l_0{-}p_0l_2: p_0l_1{-}p_1l_0).
\end{equation*}
\noindent \textit{Proof}: Most of the results can be obtained by straight forward computation. Here we just show that a point on the perpendicular bisector of $[P,Q]_-$ is equidistant from $P$ and $Q$:\\  
If $R$ is a point on this perpendicular bisector, then there exist real numbers $s$ and $t$ such that
$R^\circ = s (P - Q)^\circ + t (P^\circ \times Q^\circ).$ Then, 
\begin{equation*}
\begin{split}
|R^\circ\cdot P^\circ| = |s(P - Q)^\circ \cdot P^\circ| &= |s(1- Q^\circ \cdot P^\circ)|/\|P^\circ-Q^\circ\| \\
&= |s(1- P^\circ \cdot Q^\circ)|/\|P^\circ-Q^\circ\| = |R^\circ\cdot Q^\circ|.\quad \Box               
\end{split}
\end{equation*}\vspace*{0.5mm} 

\subsubsection{} \textbf{Parallel line through a point}\hspace*{\fill} \\
\noindent\hspace*{5mm}%
Given a line $l = (l_0:l_1:l_2)^*$ and a point $P = (p_0:p_1:p_2) \ne l^\delta$, the \textit{parallel to} $l$ \textit{through} $P$, $\textrm{par}(l,P)$, is the line $\textrm{perp}(\textrm{perp}(l,P),P)$ (cf. [27]):\\
\begin{equation*}
\begin{split}
&\textrm{par}(l,P) = \\
&(p_0 (l_1 p_1{+}l_2 p_2){-}l_0 (p_1^2{+}p_2^2): p_1 (l_0 p_0{+}l_2 p_2){-}l_1 (p_0^2{+}p_2^2): p_2 (l_0 p_0{+}l_1 p_1){-}l_2 (p_0^2{+}p_1^2))^\ast.
\end{split}\vspace*{5 mm}
\end{equation*}\vspace*{0.5mm}

\subsubsection{} \textbf{Reflections}\hspace*{\fill} \\
\noindent\hspace*{5mm}%
The mirror image of a point $P = (p_0:p_1:p_2)$ in a point $S = (s_0:s_1:s_2)$ is the point $Q = (q_0:q_1:q_2)$ with
\begin{equation*}
\begin{split}
&q_0 = {p_0( s_0^2-s_1^2-s_2^2)+2s_0(p_1s_1+p_2s_2)},\\ &q_1 = {p_1(-s_0^2+s_1^2-s_2^2) + 2s_1(p_0s_0+p_2 s_2)},\\&q_2 = {p_2(-s_0^2-s_1^2+s_2^2)+2s_2(p_0 s_0+p_1s_1). }
\end{split}
\end{equation*}

\noindent \textit{Proof}: 
By using a computer algebra system (CAS) it can be confirmed that 
\begin{equation*}
\frac{(p_0,p_1,p_2)}{\|(p_0,p_1,p_2)\|} + \frac{(q_0,q_1,q_2)}{\|(q_0,q_1,q_2)\|} = \frac {2(p_0s_0+p_1s_1+p_2s_2)}{\|(p_0,p_1,p_2)\|\|(s_0,s_1,s_2)\|} (s_0,s_1,s_2).
\end{equation*} 
From this equation results that $S$ is a midpoint of either $[P,Q]_+$ or $[P,Q]_-\;\;\Box \vspace*{3.5mm}$.

\noindent \textit{Remark:} A reflexion in a point $S$ can also be interpreted as 
\begin{itemize} 
\item[-]a rotation about $S$ through an angle of $\frac {\pi}{2},$ 
\item[-]a reflexion in the line $S^\delta.$
\end{itemize} 

\subsubsection{} \textbf{Circles}\hspace*{\fill} \\
\noindent\hspace*{5mm}%
For two points $M = (m_0:m_1,m_2)$ and $P = (p_0:p_1:p_2)$, the circle $\mathcal{C}(M,P)$ with center $M$ through the point $P$ consists of all points $X = (x_0:x_1:x_2)$ with $X^\circ{\cdot}\,M^\circ = P^\circ{\cdot}\,M^\circ.$ Thus, the coordinates of $X$ must satisfy the quadratic equation 
\begin{equation*}
(x_0^{ 2}+x_1^{ 2}+x_2^{ 2})(m_0p_0+m_1p_1+m_2p_2)^2 - (m_0x_0+m_1x_1+m_2x_2)^2(p_0^2+p_1^2+p_2^2) = 0.
\end{equation*}

\section{The use of barycentric coordinates}\hspace*{\fill} \\

\begin{figure}[!htbp]
\includegraphics[height=8cm]{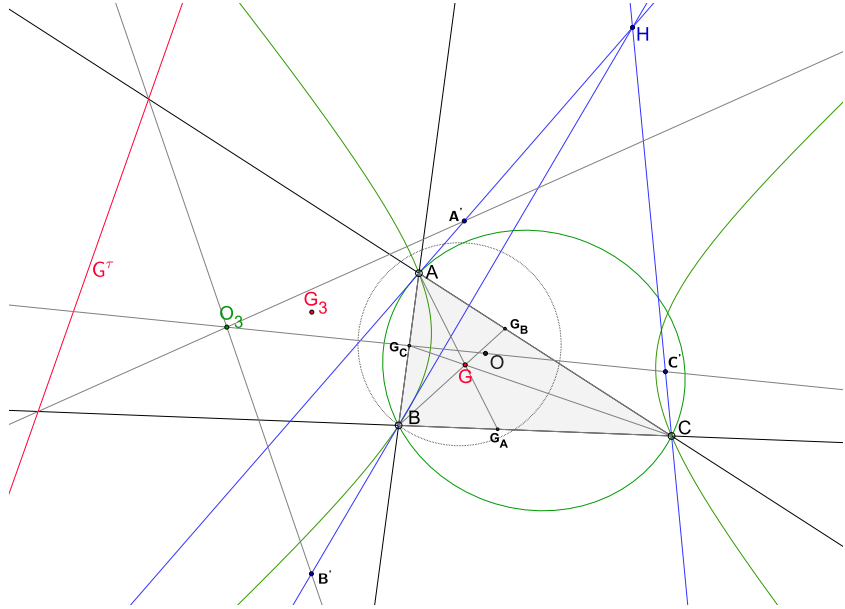}
\caption{The picture shows the triangle $\Delta_0$, the dual triple $A'B'C'$ of $ABC$ and the orthocenter $H$, furthermore the centroids $G$ and $G_3$ of $\Delta_0$ and $\Delta_3$, the tripolar line of $G$, as well as the circumcircles of $\Delta_0$ and $\Delta_3$ together with their centers $O$ and $O_3$. Since the absolute conic $\mathcal{C}$ has no real points, the dotted circle  $\tilde{\mathcal{C}} := \{(x_0:x_1:x_2) |\, x_0^{\;2} = x_1^{\;2}+x_2^{\;2}\}$ serves as a substitute for constructions. For example, the pole $A'$ of the line $B \vee C$ with respect to $\mathcal{C}$ can be obtained as follows: Construct the pole of $B \vee C$ with respect to $\tilde{\mathcal{C}}$, then its mirror image in the center $(1:0:0)$ of $\tilde{\mathcal{C}}$ 
is $A'$.\newline
The figures were created with the software program GeoGebra [29].}
\end{figure}

\subsection{} \textbf{Triangles, their sites and inner angles}\hspace*{\fill} \\
\noindent\hspace*{5mm}%
We now fix a reference triple $ABC$ of non-collinear points in $\mathcal{P}$. 
The set $\mathcal{P} - (A{\vee}B \cup B{\vee}C \cup C{\vee}A)$ consists of four connected components. Their closures are called \textit{triangles}. Thus, there are four triangles $\Delta_0, \Delta_1, \Delta_2, \Delta_3$ that share the same vertices $A, B, C$.  Inside each triangle, there is exactly one of the four points ${G_0 := A{+}B{+}C}$, ${G_1 :={-}A{+}B{+}C}$, ${G_2 := A{-}B{+}C, G_3 := A{+}B{-}C}$, and we enumerate the triangles such that $G_i \in \Delta_i$ for $i = 0,1,2,3.$ In this case, the point $G_i$ is the \textit{centroid} of the triangle $\Delta_i$ (for $i = 0,1,2,3).$\hspace*{\fill} \\
Besides the vertices, the four triangles $\Delta_i$ have all the same \textit{sidelines} $B \vee C,$ $C \vee A,$ $A \vee B$, but they do not have the same \textit{sides}. For example, the sides of $\Delta_0$ are $[B,C]_+, [C,A]_+,$ $[A,B]_+$, while $\Delta_1$ has the sides $[B,C]_+,  [C,A]_- , [A,B]_-$. The lengths of the sides of the triangle $\Delta_0$ we denote by $a_0 := \mu([B,C]_+), b_0:= \mu([C,A]_+), c_0:= \mu([A,B]_+)$. The side lengths of the other triangles $\Delta_i$ are named accordingly; for example, $a_1 := \mu([B,C]_+), b_1:= \mu([C,A]_-), c_1:= \mu([A,B]_-)$. \\
We introduced angles as a subset of the pencil of lines through a point. The (inner) angles of $\Delta_0$ are $\angle_{+} BAC,  \angle_{+} CBA,\angle_{+} ACB$, the angles of $\Delta_1$ are $\angle_{+} BAC$,  $\angle_{-} CBA$, $\angle_{-} ACB,$ etc.
The measures of these angles are $\alpha_0 = \mu(\angle_{+} BAC)$, $\beta_0 = \mu(\angle_{+} CBA)$, $\gamma_0 = \mu(\angle_{+} ACB)$, $\alpha_1 = \mu(\angle_{+} BAC)$, $\beta_1 = \mu(\angle_{-} CBA) = \pi - \mu(\angle_{+} CBA)$, etc. \\

In the following we concentrate mainly on the triangle $\Delta_0$. After having calculated the coordinates of triangle centers for this triangle, the results can be easily transferred to the triangles $\Delta_i, i > 0$. To simplify the notation, we write $a, b, c, \alpha, \beta,\gamma$ instead of $a_0, b_0, c_0, \alpha_0, \beta_0,\gamma_0.$\\
\noindent To shorten formulas, we will use abbreviations:\\
For $x \in \mathbb{R}$ define $\textrm{c}_x := \cos x$ and $\textrm{s}_x := \sin x$.\\
The semiperimeter of the triangle $\Delta_0$ is $s := (a+b+c)/2$.\\
We put \vspace*{-5mm}
\begin{equation*}
\begin{split}
\textrm{S}_A &:=  \textrm{c}_a -\textrm{c}_b \textrm{c}_c = \textrm{s}_{s} \textrm{s}_{s-a} - \textrm{s}_{s-b} \textrm{s}_{s-c},  \\
\textrm{S}_B &:= \textrm{c}_b - \textrm{c}_c \textrm{c}_a   ,\;\;\textrm{S}_C :=  \textrm{c}_c - \textrm{c}_a \textrm{c}_b,  
\end{split}
\end{equation*}
and for the barycentric coordinates of a point $P$ with respect to the reference triple $\Delta = ABC$ we use the short form $[p_1:p_2:p_3]$ instead of $[p_1:p_2:p_3]_\Delta.$

\subsection{} \textbf{Rules of elliptic trigonometry} \hspace*{\fill} \vspace*{-3mm}\\

In the following, we will make use of the following rules of elliptic trigonometry:\\

\noindent Cosine rules: 
\begin{equation*}
\begin{split}
&\cos \alpha = \frac {\cos a - \cos b \,\cos c} {\sin b \,\sin c}\, ,\, \cos a = \frac {\cos \alpha + \cos \beta \,\cos \gamma} {\sin \beta\sin \gamma}
 = 1 +  \frac{2 \sin \epsilon \, \sin \epsilon\text{-}\alpha}{\sin \beta \, \sin \gamma}, \\
&\text{where}\; 2\epsilon = \alpha{+}\beta{+}\gamma{-}\pi \;\text{is the \textit{excess} of\;} \Delta_0.
\end{split}
\end{equation*}
Sine rule: 
\begin{equation*}
\hspace*{-16mm}\frac {\sin \alpha} {\sin a} = \frac {\sin \beta} {\sin b} = \frac {\sin \gamma} {\sin c} = \frac {|\textrm{S}|}{\sin a \,\sin b \,\sin c}, 
\end{equation*}
with \hspace*{22mm}$\textrm{S} = \textrm{S}(\Delta) = det (\begin{pmatrix} A^\circ \\ B^\circ\\ C^\circ\end{pmatrix})$\\
and 
\begin{equation*}
\begin{split}
\hspace*{11mm}|\textrm{S}| &= {\sqrt{1 - \cos^2 a - \cos^2 b- \cos^2 c + 2\cos a \,\cos b\, \cos c}}\\
&= 2 \sqrt{\sin s\, \sin (s-a)\, \sin (s-b) \,\sin (s-c)}.
\end{split}
\end{equation*}

\noindent \textit{Proof}: 
	\[ 
\begin{split}
	\cos \alpha &= \frac {(A^\circ\times B^\circ){\cdot}(A^\circ\times C^\circ)}{\|A^\circ\times B^\circ\| \|A^\circ\times C^\circ\|} = \frac { (A^\circ{\cdot} A^\circ)(B^\circ{\cdot}\,C^\circ) - (A^\circ{\cdot}\,C^\circ)(B^\circ{\cdot} A^\circ)}{\textrm{s}_c \textrm{s}_b  }\\
	&= \frac {\textrm{c}_a - \textrm{c}_b \textrm{c}_c} {\textrm{s}_b \textrm{s}_c  }.
\end{split} 
\]
Before giving a proof for the second cosine rule, we prove the sine rule.	
\begin{equation*}
\begin{split}
\frac {\sin \alpha} {\sin a} &= \frac {1} {\textrm{s}_a} \frac {\|(A^\circ\times B^\circ)\times (A^\circ\times C^\circ)\|}{\|A^\circ\times B^\circ\| \|A^\circ\times C^\circ\|} = \frac {| \textrm{S}(\Delta)| }{\textrm{s}_a \,\textrm{s}_b \,\textrm{s}_c }\\
\sin^2 \alpha &= 1 - \cos^2 \alpha = 1 - \big( \frac {\textrm{c}_a - \textrm{c}_b \textrm{c}_c} {\textrm{s}_b \textrm{s}_c  } \big )^2 
= \frac {{1 - \textrm{c}_a^2 - \textrm{c}_b^2- \textrm{c}_c^2 + 2\textrm{c}_a \textrm{c}_b \textrm{c}_c}} {\textrm{s}_b^2 \textrm{s}_c^2}\\
\sin^2 \alpha &= (1+\cos \alpha)(1-\cos \alpha) =(1 +  \frac {\textrm{c}_a - \textrm{c}_b \textrm{c}_c} {\textrm{s}_b \textrm{s}_c  } )(1 -  \frac {\textrm{c}_a - \textrm{c}_b \textrm{c}_c} {\textrm{s}_b \textrm{s}_c } )\\
&= \frac{\textrm{c}_a - \textrm{c}_{b+c}}{\textrm{s}_b \textrm{s}_c} \,\frac{\textrm{c}_{b-c} - \textrm{c}_{a}}{\textrm{s}_b \textrm{s}_c} 
= \frac{4 \textrm{s}_{s} \textrm{s}_{s-a}\textrm{s}_{s-b} \textrm{s}_{s-c}}{\textrm{s}_b^{\,2} \textrm{s}_c^{\,2}}.
\end{split} 
\]
Proof of the second cosine rule:
	\[
\begin{split}	
	\frac {\cos \alpha + \cos \beta \,\cos \gamma} {\sin \beta\sin \gamma} &=\frac {\frac {\cos a - \cos b \,\cos c} {\sin b \,\sin c} + \frac {\cos b - \cos c \,\cos a} {\sin c \,\sin a} \,{\frac {\cos c - \cos a \,\cos b} {\sin a \,\sin b}}} {\sin \beta \,\sin \gamma}\\
&= 	\frac {\textrm{c}_a (1 - \textrm{c}_a^2 - \textrm{c}_b^2- \textrm{c}_c^2 + 2\textrm{c}_a \textrm{c}_b \textrm{c}_c)} {\textrm{s}_\beta \,\textrm{s}_\gamma \textrm{s}_a^{\;2}  \textrm{s}_b \textrm{s}_c} = \textrm{c}_a 
\end{split} 	
	\]\vspace*{-1mm}
	\[
\begin{split}	\cos a &= 1 + \frac{\cos \alpha + \cos \beta \cos \gamma - \sin \beta \sin \gamma}{\sin \beta \sin \gamma}\\
&= 1 + \frac{\cos \alpha + \cos \beta{+}\gamma}{\sin \beta \sin \gamma}\\
&= 1 + \frac{ 2\,\cos \frac 1 2(\alpha{+}\beta{+}\gamma)\;\cos \frac 1 2 (\beta{+}\gamma{-}\alpha)}{\sin \beta \sin \gamma} 
= 1 +  \frac{2 \sin \epsilon \, \sin \epsilon\text{-}\alpha}{\sin \beta  \sin \gamma} \;\;\; \;\;\;\;\;\;\Box	
\end{split} 	
	\]
	
A detailed collection of spherical trigonometry formulas, including their proofs, can be found in [4, 24].
	
\subsection{Calculations with barycentric coordinates} \hspace*{\fill} \\

\noindent By using barycentric coordinates, many concepts can be transferred directly from metric affine (for example euclidean) geometry to elliptic geometry. \hspace*{\fill} 

\subsubsection{} \textbf{The distance between two points in barycentric coordinates} \hspace*{\fill} \\
\noindent\hspace*{5mm}We introduce the matrix 
\begin{equation*}
\mathfrak{T} = (\mathfrak{t}_{ij})_{i,j=1,2,3}:= \begin{pmatrix} A^\circ\cdot A^\circ &A^\circ\!\cdot B^\circ &A^\circ\!\cdot C^\circ\\B^\circ\!\cdot A^\circ &B^\circ\!\cdot B^\circ &B^\circ \!\cdot C^\circ \\C^\circ\!\cdot A^\circ &C^\circ\!\cdot B^\circ &C^\circ\!\cdot C^\circ \end{pmatrix} = \begin{pmatrix}1 &\cos c &\cos b\\ \cos c & 1 & \cos a\\ \cos b & \cos a & 1 \end{pmatrix},
\end{equation*}
which we call the \textit{characteristic matrix} of $\Delta$. \\
Besides the dot product $\cdot$ for 3-vectors we introduce another scalar product $\ast$:
\begin{equation*}
\begin{split}
  &(p_1,p_2,p_3) \ast (q_1,q_2,q_3) \\
	&= (p_1,p_2,p_3)\, \mathfrak{T}\!\begin{pmatrix} q_1 \\ q_2\\ q_3 \end{pmatrix}\\
 &= p_1q_1+p_2q_2+p_3q_3+(p_2q_3{+}p_3q_2)\cos a +(p_3q_1{+}p_1q_3)\cos b +(p_1q_2{+}p_2q_2)\cos c
\end{split}
\end{equation*}
and use the abbreviations
\begin{equation*}
\begin{split}
&(p_1,p_2,p_3)^2 := (p_1,p_2,p_3)\cdot(p_1,p_2,p_3), (p_1,p_2,p_3)^{\ast2} := (p_1,p_2,p_3)\ast(p_1,p_2,p_3),\\
&\|(p_1,p_2,p_3)\|_\ast:= \sqrt{(p_1,p_2,p_3)^{\ast2}\,}.
\end{split}
\end{equation*}

With the help of these products and the resulting norms, the distance between two points $P = {[p_1:p_2:p_3]}$ and $Q = [q_1:q_2:q_3]$ can be calculated as follows:
\begin{equation*}
\begin{split}
d(P,Q) &= \arccos\;\frac {|(p_1A^\circ+p_2B^\circ+p_3C^\circ)\cdot(q_1A^\circ+q_2B^\circ+q_3C^\circ)|}{{\|p_1A^\circ+p_2B^\circ+p_3C^\circ\|}\;\;{\|q_1A^\circ+q_2B^\circ+q_3C^\circ\|}}\\
       &= \arccos \,\frac{|(p_1,p_2,p_3) \ast (q_1,q_2,q_3)|}{\;\,\|{(p_1,p_2,p_3)\|_\ast\,\|(q_1,q_2,q_3)\|_\ast}}.
\end{split}
\end{equation*}

\subsubsection{} \textbf{Circles} \hspace*{\fill} \\
\noindent\hspace*{5mm}For two points $M = [m_1:m_2:m_3]$ and $P = [p_1:p_2:p_3]$, the circle $\mathcal{C}(M,P)$ with center $M$ through the point $P$ consists of all points $X = [x_1:x_2:x_3]$ with 
\begin{equation*}
((m_1,m_2,m_3){\ast}(x_1,x_2,x_3))^2(p_1,p_2,p_3)^{\ast2} = ((m_1,m_2,m_3){\ast}(p_1,p_2,p_3))^2(x_1,x_2,x_3)^{\ast2}.
\end{equation*}

\subsubsection{} \textbf{Lines} \hspace*{\fill} \\
\noindent\hspace*{5mm}Given two different points $P = [p_1\!:\!p_2\!:\!p_3]$ and $Q = [q_1\!:\!q_2\!:\!q_3]$, a third point $X = [x_1\!:\!x_2\!:\!x_3]$ is a point on $P \vee Q$ exactly when 
\begin{equation*}
((p_1,p_2,p_3)\times(q_1,q_2,q_3))\cdot(x_1,x_2,x_3) = 0.
\end{equation*}
\noindent\hspace*{5mm}If $R \vee S$ is a line through $R = [r_1\!:\!r_2\!:\!r_3]$ and $S = [s_1\!:\!s_2\!:\!s_3]$, different from $P \vee Q$, then both lines meet at a point 
$T = [t_1\!:\!t_2\!:\!t_3]$ with 
\begin{equation*}
(t_1,t_2,t_3) = ((p_1,p_2,p_3)\times ((q_1,q_2,q_3))\times((r_1,r_2,r_3){\times}(s_1,s_2,s_3)).
\end{equation*}

\subsubsection{} \textbf{The midpoint of a segment} \hspace*{\fill} \\
\noindent\hspace*{5mm}Let  $P = [p_1:p_2:p_3]$ and $Q = [q_1:q_2:q_3]$ are two different points. We can assume that the triples $ \boldsymbol{p} =  (p_1,p_2,p_3)$ and $ \boldsymbol{q} =  (q_1,q_2,q_3)$ are both strictly positive with respect to the lexicographic order.
Then the midpoints of the two segments $[P,Q]_\pm$ are \vspace*{-2mm}
\begin{equation*}
[\frac {p_1}{\|\boldsymbol{p}\|_\ast}\pm \frac {q_1}{\|\boldsymbol{q}\|_\ast}:\frac {p_2}{\|\boldsymbol{p}\|_\ast} \pm \frac {q_2}{\|\boldsymbol{q}\|_\ast}:\frac {p_3}{\|\boldsymbol{p}\|_\ast} \pm \frac {q_3}{\|\boldsymbol{q}\|_\ast}].
\end{equation*}

\subsubsection{} \textbf{The dual $\Delta'$ of $\Delta$}\hspace*{\fill} \\
\noindent\hspace*{5mm}We put $A' := (B{\vee}C)^\delta, B' := (C{\vee}A)^\delta, C' := (A{\vee}B)^\delta$. The triple $\Delta' = A'B'C'$ is called the dual of $\Delta$.
The points $A', B', C'$ can be written in terms of barycentric coordinates as follows:
\begin{equation*}
\begin{split}
A' &= [1- \textrm{c}_a^2 :\textrm{c}_a \textrm{c}_b - \textrm{c}_c : \textrm{c}_c \textrm{c}_a - \textrm{c}_b] = [\textrm{s}_a^2: -\textrm{S}_C:-\textrm{S}_B]\\
B' &= [\textrm{c}_a \textrm{c}_b - \textrm{c}_c :1- \textrm{c}_b^2 : \textrm{c}_b \textrm{c}_c - \textrm{c}_a] = [-\textrm{S}_C:\textrm{s}_b^2: -\textrm{S}_A]\\
C' &= [\textrm{c}_c \textrm{c}_a - \textrm{c}_b : \textrm{c}_b \textrm{c}_c - \textrm{c}_a : 1- \textrm{c}_c^2] = [-\textrm{S}_B: -\textrm{S}_A: \textrm{s}_c^2]
\end{split}
\end{equation*}
\noindent \textit{Proof}: Up to a real factor $1/\textrm{S}$, the characteristic matrix $\mathfrak{T}$ is the matrix that transforms $\begin{pmatrix} B^\circ{\times}C^\circ\\C^\circ{\times}A^\circ\\A^\circ{\times}B^\circ \end{pmatrix}$ onto $\begin{pmatrix} A^\circ\\B^\circ\\C^\circ \end{pmatrix}$:
\begin{equation*}
\begin{pmatrix} A^\circ\\B^\circ\\C^\circ \end{pmatrix} =  \frac 1 {\textrm{S}(\Delta)} \mathfrak{T} \begin{pmatrix} B^\circ{\times}C^\circ\\C^\circ{\times}A^\circ\\A^\circ{\times}B^\circ \end{pmatrix}.
\end{equation*}

\noindent The matrix $\mathfrak{T}^{-1}$ of the inverse transformation is
\begin{equation*}
\begin{split}
\mathfrak{T}^{-1} &= \frac {1}{\textrm{S}} \begin{pmatrix} \mathfrak{t}_{22}\mathfrak{t}_{33}-\mathfrak{t}_{23}^2&\mathfrak{t}_{23}\mathfrak{t}_{31}-\mathfrak{t}_{12}\mathfrak{t}_{33}&{\mathfrak{t}_{12}\mathfrak{t}_{23}-\mathfrak{t}_{31}\mathfrak{t}_{22}}\\ 
\mathfrak{t}_{23}\mathfrak{t}_{31}-\mathfrak{t}_{12}\mathfrak{t}_{33}&\mathfrak{t}_{33}\mathfrak{t}_{11}-\mathfrak{t}_{31}^2&\mathfrak{t}_{31}\mathfrak{t}_{12}-\mathfrak{t}_{13}\mathfrak{t}_{11} \\
\mathfrak{t}_{12} \mathfrak{t}_{23}-\mathfrak{t}_{31}\mathfrak{t}_{22}&\mathfrak{t}_{31}\mathfrak{t}_{12}-\mathfrak{t}_{13}\mathfrak{t}_{11}&\mathfrak{t}_{11}\mathfrak{t}_{22}-\mathfrak{t}_{12}^2 \end{pmatrix}\\
&= \frac {1}{\textrm{S}} \begin{pmatrix} 1-\textrm{c}_a^2&\textrm{c}_a\textrm{c}_b-\textrm{c}_c&{\textrm{c}_c\textrm{c}_a-\textrm{c}_b}\\ 
\textrm{c}_a\textrm{c}_b-\textrm{c}_c&1-\textrm{c}_b^2&\textrm{c}_b\textrm{c}_c-\textrm{c}_a\\
\textrm{c}_c \textrm{c}_a-\textrm{c}_b&\textrm{c}_b\textrm{c}_c-\textrm{c}_a&1-\textrm{c}_c^2 \end{pmatrix},
\end{split}
\end{equation*}  \\ 
which proves the statement. $\;\;\Box$

\subsubsection{} \textbf{The dual of a point and the dual of a line}\hspace*{\fill} \\
\noindent\hspace*{5mm}The dual line $P^\delta $ of a point $P = [p_1:p_2:p_3]$ has the equation (in barycentric coordinates)\vspace*{-4mm}
	\[(p_1,p_2,p_3) \mathfrak{T}\!\begin{pmatrix} x_1 \\ x_2\\ x_3 \end{pmatrix} = 0.\vspace*{1mm}
\]
If $\ell$ is a line with equation ${\ell_1 x_1 + \ell_2 x_2 + \ell_3 x_3 = 0}$,
then its dual is the point $R = [r_1:r_2:r_3]$ with 
$(r_1,r_2,r_3)=(\ell_1,\ell_2,\ell_3)\mathfrak{T}^{-1}.$

\subsubsection{} \textbf{The angle bisectors of two lines}\hspace*{\fill} \\
\noindent\hspace*{5mm}Let $k: k_1 x_1 + k_1 x_1 + k_1 x_1 = 0$ and $\ell: \ell_1 x_1 + \ell_1 x_1 + \ell_1 x_1 = 0$ be two different lines, then their two angle bisectors are
	$m: m_1 x_1 + m_1 x_1 + m_1 x_1 = 0$ and $n: n_1 x_1 + n_1 x_1 + n_1 x_1 = 0,$
with
	\[ 
\begin{split}	
	&(m_1,m_2,m_3) = \sqrt{(\boldsymbol{\ell}\mathfrak{T}^{-1})\cdot \boldsymbol{\ell}}\;\boldsymbol{k} + \sqrt{(\boldsymbol{k}\mathfrak{T}^{-1})\cdot \boldsymbol{k}}\;\boldsymbol{\ell},\\ 	
	&(n_1,n_2,n_3) = \sqrt{(\boldsymbol{\ell}\mathfrak{T}^{-1})\cdot \boldsymbol{\ell}}\;\boldsymbol{k} - \sqrt{(\boldsymbol{k}\mathfrak{T}^{-1})\cdot \boldsymbol{k}}\;\boldsymbol{\ell},\\
	&\boldsymbol{k} = (k_1,k_2,k_3), \boldsymbol{\ell} = (\ell_1,\ell_2,\ell_3).
\end{split}
\]

\subsubsection{} \textbf{Chasles' Theorem}\hspace*{\fill} \\
If the dual $\Delta' = A'B'C'$ of $\Delta$ is different from $\Delta$, then $\Delta'$ and $\Delta$ are a perspective; the center of perspective is \begin{equation*}
\begin{split}
H :=& [\frac{1}{\textrm{c}_b \textrm{c}_c- \textrm{c}_a}:\frac{1}{\textrm{c}_b \textrm{c}_a-\textrm{c}_b}:\frac{1}{\textrm{c}_a \textrm{c}_b-\textrm{c}_c}]\\
=& [\frac{1}{\textrm{S}_A}:\frac{1}{\textrm{S}_B}:\frac{1}{\textrm{S}_C}]\;.
\end{split}
\end{equation*}\\
\noindent \textit{Proof:} We prove the following three statements:\\
If $A'\ne A, B'\ne B, C'\ne C$, then the lines $A \vee A', B \vee B', C \vee C'$ are concurrent at $H$.
If $A = A'$, but $B'\ne B$ and $C'\ne C$, the lines $B \vee B', C \vee C'$ meet at $H$.
If $A=A'$ and $B=B'$, then $C=C'$.\\
Using 2.3.3, the first two theses statements can be verified by calculation.\\
$A=A'$ is equivalent with $d(A,B) = d(A,C) = \pi$ and $B=B'$ equivalent with $d(B,C) = d(B,A) = \pi$. Thus, if $A=A'$ and $B=B'$, then $d(C,A) = d(C,B) = \pi$ and therefore $C = C'$.\\
\noindent\hspace*{5mm}The point $H$ is called the \textit{orthocenter} of $\Delta$.\\
\noindent\hspace*{5mm}We define: A triple $PQR$ of three points $P,Q,R$ is a \textit{perspective triple} with perspector $S$ if the triples $ABC$ and $PQR$ are perspective and $S$ is the perspective center.

\subsubsection{} \textbf{Pedals and antipedals of a point}\hspace*{\fill} \\
\noindent\hspace*{5mm}The \textit{pedals} of a point $P = [p_1:p_2:p_3]$ on the sidelines of $\Delta$ are notated $A_{[P]}, B_{[P]}, C_{[P]}$. We calculate the barycentric coordinates $[q_1:q_2:q_3]$ of $A_{[P]} $:
\begin{equation*} 
\begin{split}
(q_1,q_2,q_3) &= ((p_1,p_2,p_3) \times (\textrm{s}_a^2,-\textrm{S}_C ,-\textrm{S}_B))\times ((0,1,0) \times (0,0,1)) \\
&= (0, p_1\textrm{S}_C + p_2\textrm{s}_a^2, p_1\textrm{S}_B + p_3\textrm{s}_a^2).
\end{split}
\end{equation*}
Similarly the coordinates of the other two pedals can be calculated:
	\[
\begin{split}	B_{[P]} &= [p_2 \textrm{S}_C + p_1 \textrm{s}_b^{\;2} : 0 : p_2 \textrm{S}_A + p_3 \textrm{s}_b^{\;2}]\\
              C_{[P]} &= [p_3 \textrm{S}_B + p_1 \textrm{s}_c^{\;2} : p_3 \textrm{S}_A + p_2 \textrm{s}_c^{\;2} : 0].
\end{split}\vspace*{1mm}  
\]
\noindent\hspace*{5mm}Define the \textit{antipedal points} $A^{[P]}, B^{[P]}, C^{[P]}$ of $P$ by\vspace*{1mm}\\
\centerline{$A^{[P]} := \textrm{perp}(B{\vee}P,B) \wedge  \textrm{perp}(C{\vee}P,C)$ and  $B^{[P]}, C^{[P]}$ cyclically.}\vspace*{1mm}\\   
Straightforward calculation gives 
\begin{equation*} 
 A^{[P]} = [-1:\frac{p_2 \textrm{S}_C + p_1 \textrm{s}_b^2}{p_1 \textrm{S}_C  + p_2 \textrm{s}_a^2}:\frac{p_3 \textrm{S}_B + p_1 \textrm{s}_c^2}{p_1 \textrm{S}_B + p_3 \textrm{s}_a^2}].\vspace*{1mm}
\end{equation*} 
A special case: For $P = H$ we get 
	\[A^{[H]} = [-\textrm{c}_a:\textrm{c}_b:\textrm{c}_c],\; B^{[H]} = [\textrm{c}_a:-\textrm{c}_b:\textrm{c}_c], \;A^{[H]} = [\textrm{c}_a:\textrm{c}_b:-\textrm{c}_c]. 
	\]
	
\subsubsection{} \textbf{Cevian and anticevian triangles}\hspace*{\fill} 
\noindent\hspace*{5mm}%

If $P = {[p_1:p_2:p_3]}$ is a point different from $A,B,C$, then the lines $P\vee A,$ $P\vee B,$ $P\vee C$ are called the \textit{cevians} of $P$. The cevians meet the sidelines $a, b, c$ in $A_{P} :={[0:p_2:p_3]}, B_{P} := {[p_1:0:p_3]}, C_{P} := {[p_1:p_2:0]}$, respectively. These points are called the \textit{traces} of $P$.  
The points $A^{P} := {[-p_1:p_2:p_3]}, B^{P} := {[p_1:-p_2:p_3]}, C^{P} := {[p_1:p_2:-p_3]}$ are called \textit{harmonic associates} of $P$.  \\

We now assume that $P$ is not a point on a sideline of $ABC$ and define: The \textit{cevian triangle} of $P$ with respect to $\Delta_0$ is the triangle with vertices $A_P, B_P, C_P$ which contains the point $[|p_1|:|p_1|:|p_1|]$, the cevian triangle of $P$ with respect to $\Delta_1$ is the triangle with vertices $A_P, B_P, C_P$ which contains the point $[-|p_1|:|p_1|:|p_1|]$, and so on.
Furthermore, we define: The \textit{anticevian triangle} of $P$ with respect to $\Delta_0$ is the triangle with vertices $A^P, B^P, C^P$ that has all points $A,B,C$ on its sides. The anticevian triangle of $P$ with respect to $\Delta_1$ has the same vertices, but only the point $A$ is on one of its sides, while the points $B$ and $C$ are not.

\noindent A special case: \\
The traces of $G_i$ are the midpoints of the sides of $\Delta_i$, the cevians of $G_i$ are (therefore) called the \textit{medians} of $\Delta_i$. $G_i$ itself is called the \textit{centroid} of $\Delta_i$, and the cevian triangle of $G_i$ with respect to $\Delta_i$ is called the \textit{medial triangle} of $\Delta_i$. \\
$\Delta_0$ is, in general, not the medial triangle of the anticevian triangle of $G_0$. (The same applies to the other triangles $\Delta_i.$) The \textit{antimedial triangle} of $\Delta_0$ is the anticevian triangle with respect to $\Delta_0$ of the point $G^+ :=[\cos a : \cos b : \cos c]$. In the last subsection it was shown that the vertices of this triangle also form the antipedal triple of $H$. We now prove that the anticevian triangle with respect to $\Delta_0$ of the point $G^+$ has $G^+$ as its centroid, cf. Wildberger [27, 28]. \vspace*{1mm}\\
\textit{Proof:}  We know already that the points $A, B, C$ lie on the sides of the anticevian triangle of $G^+$ with respect to $\Delta_0$. Now we show that $A$ is equidistant from $B^{G^+}$ and $C^{G^+}$ by proving the equation $A = B^{G^+}\!{+}\,C^{G^+}$: Define the vectors $\boldsymbol{p}$ and $\boldsymbol{q}$ by $\boldsymbol{p} := (\textrm{c}_a, -\textrm{c}_b, \textrm{c}_c)$ and $\boldsymbol{q} := (\textrm{c}_a, \textrm{c}_b, -\textrm{c}_c)$. Since $\boldsymbol{p}^{\ast2} = \boldsymbol{q}^{\ast2}$, we have $\boldsymbol{p}/\|\boldsymbol{p}\|_\ast + \boldsymbol{q}/\|\boldsymbol{q}\|_\ast = (2\textrm{c}_c ,0,0)/\|\boldsymbol{p}\|_\ast.$ \\
In the same way it is shown that the points $B$ and $C$ are the midpoints of the corresponding sides of the anticevian triangle of $G^+$ with respect to $\Delta_0.\;\;\; \Box$\vspace*{1mm}\\
\textit{Remark:} Wildberger's name for the antimedial triangle is \textit{double triangle} [27, 28].

\subsubsection{} \textbf{Tripolar and tripole}\hspace*{\fill} \\
\noindent\hspace*{5mm}
Given a point $P = [p_1:p_2:p_3]$ different from $A, B, C$, then the point $[0:-p_2:p_3]$ is the harmonic conjugate of $A_P$ with respect to $\{B,C\}$. Correspondingly, the harmonic conjugates of the traces of $P$ on the other sidelines are $[-p_1:0:p_3]$  and $[p_1:-p_2:0]$. These three harmonic conjugates are collinear; the equation of the line $l$ is 
	\[ p_2p_3 x_1 + p_3p_1 x_2 + p_1p_2x_3 = 0.
\]
This line is called the \textit{tripolar line} or the \textit{tripolar} of $P$ and we denote it by $P^\tau$. $P$ is the \textit{tripole} of $l$ and we write $P = l^\tau$.\\

We calculate the coordinates of the dual point of the tripolar of $P$ and get 
\begin{equation*} 
P^{\tau \delta} = [p_1 (p_2 \textrm{S}_B + p_3 \textrm{S}_C) - p_2 p_3 \textrm{s}_a^2 
: p_2 (p_3 \textrm{S}_C + p_1 \textrm{S}_A)- p_3 p_1 \textrm{s}_b^2 
: p_3 (p_1 \textrm{S}_A + p_2 \textrm{S}_B) - p_1 p_2 \textrm{s}_c^2].
\end{equation*}

\noindent\textit{Two Examples:}\\ 
\noindent\hspace*{5mm}The line with the equation
	\[\textrm{S}_A x_1 + \textrm{S}_B x_2 + \textrm{S}_C x_3 = 0,
\]
is called \textit{orthic axis}. If $H$ is different from $A, B, C$, this line is the tripolar $H^\tau$.
The dual of the ortic axis is the point
	\[H^{\star} := [2 \textrm{S}_B \textrm{S}_C - \textrm{S}_A \textrm{s}_a^2:2 \textrm{S}_C \textrm{S}_A - \textrm{S}_B \textrm{s}_b^2:2 \textrm{S}_A \textrm{S}_B - \textrm{S}_C \textrm{s}_c^2].
\]
Wildberger [24, 25] names this point \textit{orthostar}, we adopt this terminology.\\

The tripolar of $G$ has the equation $x_1 + x_2 + x_3 = 0$, the point $O:=G^{\tau\delta}$ has coordinates 
\begin{equation*} 
\begin{split}
  &[(1{-}\textrm{c}_a)(1{+}\textrm{c}_a{-}\textrm{c}_b{-}\textrm{c}_c) : (1{-}\textrm{c}_b)(1{-}\textrm{c}_a{+}\textrm{c}_b{-}\textrm{c}_c) 
: (1{-}\textrm{c}_c)(1{-}\textrm{c}_a{-}\textrm{c}_b{+}\textrm{c}_c) ]\\
= &[\textrm{s}_{a/2}^2(\textrm{s}_{a/2}^2{-}\textrm{s}_{b/2}^2{-}\textrm{s}_{c/2}^2) : \textrm{s}_{b/2}^2({-}\textrm{s}_{a/2}^2{+}\textrm{s}_{b/2}^2{-}\textrm{s}_{c/2}^2) 
: \textrm{s}_{c/2}^2({-}\textrm{s}_{a/2}^2{-}\textrm{s}_{b/2}^2{+}\textrm{s}_{c/2}^2)].\\ 
\end{split}
\end{equation*}
In the next subsection we identify this point as the \textit{circumcenter} of the triangle $\Delta_0$.

\subsubsection{} \textbf{The four classical triangle centers of $\Delta_0$} \hspace*{\fill} \\
\noindent\hspace*{5mm}
We already calculated the coordinates of two classical triangle centers of $\Delta_0$, the centroid $G = G_0$ and the orthocenter $H$ :\\
\begin{equation*} 
\begin{split}
G &= [1:1:1] \\
H &= [\frac{1}{\cos a-\cos b \cos c}:\frac{1}{\cos b-\cos c \cos a}:\frac{1}{\cos c-\cos a \cos b}] 
\end{split}
\end{equation*}\\
The other two classical centers are the center $O = O_0$ of the circumcircle and the center $I = I_0$ of the incircle.\\

The point $O=G^{\tau\delta}$ is the circumcenter of $\Delta_0$;
it can be easily checked that the pedals of  $G^{\tau\delta}$ are the traces of $G$, so $G^{\tau\delta}$ is a point on all three perpendicular bisectors ${A_{G} \vee A'},\, {B_{G} \vee B'},\, {B_{G} \vee B'}$ of the triangle sides and has, therefore, the same distance from the vertices $A, B$ and $C$. 

\noindent The radius of the circumcircle is 
\begin{equation*} 
\begin{split}
 \textrm{R} &= d(O,A)           \\
&= \arccos\, \sqrt{\Big|{\frac{\textrm{c}_a^{\;2} + \textrm{c}_b^{\;2} + \textrm{c}_c^{\;2} - 2 \textrm{c}_a \textrm{c}_b \textrm{c}_c - 1} {\textrm{c}_a^{\;2}   + \textrm{c}_b^{\;2}   + \textrm{c}_c^{\;2}  - 2 \textrm{c}_b \textrm{c}_c - 2 \textrm{c}_c \textrm{c}_a - 2 \textrm{c}_a \textrm{c}_b + 2 \textrm{c}_a +  2 \textrm{c}_b + 2\textrm{c}_c - 3}}\Big|}\;\;\\
&=  \arccos\, \sqrt{\Big|\frac{\textrm{s}_a^{\;2}{+}\textrm{s}_b^{\;2}{+}\textrm{s}_c^{\;2}-( \textrm{s}_s^{\;2}{+}\textrm{s}_{s-a}^{\;2}{+}\textrm{s}_{s-b}^{\;2}{+}\textrm{s}_{s-c}^{\;2} )} {\textrm{s}_a^{\;2}{+}\textrm{s}_b^{\;2}{+}\textrm{s}_c^{\;2}+2(\textrm{s}_{s{-}b} \textrm{s}_{s{-}c}{+}\textrm{s}_{s{-}c} \textrm{s}_{s{-}a}{+}\textrm{s}_{s{-}a} \textrm{s}_{s{-}b}{-}\textrm{s}_{s} (\textrm{s}_{s{-}a}{+}\textrm{s}_{s{-}b}{+}\textrm{s}_{s{-}c})) }\Big|}\;.
\end{split}
\end{equation*}

\noindent The equation of the circumcircle: A point $X = [x_1:x_2:x_3]$ is a point on the circumcircle of $\Delta_0$ precisely when its coordinates satisfy the equation
\begin{equation*} 
(1-\cos a) x_2 x_3 + (1-\cos b) x_3 x_1 + (1-\cos c) x_1 x_2 = 0,
\end{equation*}
which is equivalent to
\begin{equation*} 
\sin^2\!\tfrac a 2\; x_2 x_3 + \sin^2\! \tfrac b 2\; x_3 x_1 + \sin^2\!\tfrac c 2\; x_1 x_2 = 0.
\end{equation*}\\
\noindent\hspace*{5mm}The incenter $I = [\sin a: \sin b:\sin c]$ of $\Delta_0$ is the meet of the three bisectors of the (inner) angles of $\Delta_0$.\vspace*{1mm}\\
\textit{Proof:} 
We show that $I = [\sin a: \sin b:\sin c]$ is a point on the bisector of $\alpha$ by showing that the dual point $W_A$ of this bisector is
orthogonal to $I$.\\ 
We introduce the vectors $\boldsymbol{i} := \|B^\circ\times C^\circ\|A^\circ + \|C^\circ \times A^\circ\|B^\circ + \|A^\circ \times B^\circ\|C^\circ$ and $\boldsymbol{w}_A := \frac{A^\circ \times B^\circ}{\|(A^\circ \times B^\circ)\|} + \frac {A^\circ \times C^\circ}{\|(A^\circ \times C^\circ)\|}$; $\boldsymbol{i}$ is a multiple of $I^\circ$, $\boldsymbol{w}_A$ a multiple of $(W_A)^{\circ}$.
It can be easily checked that $\boldsymbol{i}\cdot \boldsymbol{w}_A = 0$.\\
In a similar way it can be proved that $I$ is a point on the other two angle bisectors.$\;\Box$  \vspace*{1mm}\\
The pedals of $I$ are 
\begin{equation*} 
\begin{split}
A_{[I]} &= [0: \cos a \,\cos b  - \sin a \,\sin b - \cos c: \cos c \,\cos a - \sin c \,\sin a - \cos b]\\  
B_{[I]} &= [\cos a \,\cos b  - \sin a \,\sin b - \cos c:0: \cos b \,\cos c - \sin b \,\sin c - \cos a]\\  
C_{[I]} &= [\cos c \,\cos a - \sin c \,\sin c - \cos b: \cos b \,\cos c - \sin b \,\sin c - \cos a : 0]
\end{split}
\end{equation*}

\noindent These three pedal points are also the traces of the \textit{Gergonne point} 
\begin{equation*} 
\begin{split}
Ge &= [\frac {1}{\cos (b{+}c) - \cos a}: \frac {1}{\cos (c{+}a) - \cos b}:\frac {1}{\cos(a{+}b) - \cos c}]\\
&=  [\frac {1}{\sin (s{-}a)}: \frac {1}{\sin (s{-}b)}:\frac {1}{\sin (s{-}c)}].  
\end{split}
\end{equation*}
The cevian triangle of $Ge$ is called the \textit{tangent triangle} of $\Delta_0$.\\

\noindent The radius of the inner circle is 
\begin{equation*} 
\textrm{r} = \arccos \;{\frac{2}{\kappa}\sin s}  \\ 
\end{equation*}
with $\kappa = \|(\textrm{s}_a,\textrm{s}_b ,\textrm{s}_c)\|_\ast =\sqrt{\textrm{s}_a^2 {+}\,\textrm{s}_b^2\, {+}\,\textrm{s}_a^2\, {+}\,2(\textrm{c}_a \textrm{s}_b \textrm{s}_c\, {+}\,\textrm{c}_b \textrm{s}_c \textrm{s}_a \,{+}\,\textrm{c}_c \textrm{s}_a \textrm{s}_b)}\;.$\\

\noindent \textit{Proof:} 
\begin{equation*} 
\begin{split}
\cos d(I,A_{[I]}) &= \sqrt{ \frac {((\textrm{s}_a,\textrm{s}_b,\textrm{s}_c) \ast (0,\textrm{c}_a \textrm{c}_b - \textrm{s}_a \textrm{s}_b - \textrm{c}_c, \textrm{c}_c \textrm{c}_a - \textrm{s}_c \textrm{s}_a - \textrm{c}_b))^2}{{(\textrm{s}_a,\textrm{s}_b ,\textrm{s}_c)^{{\ast}2} (0,\textrm{c}_a \textrm{c}_b - \textrm{s}_a \textrm{s}_b - \textrm{c}_c, \textrm{c}_c \textrm{c}_a - \textrm{s}_c \textrm{s}_a - \textrm{c}_b)^{\ast 2} }}}\\ 
&= \sqrt{ \frac { (2\,\textrm{s}_a (1-(\textrm{c}_a (\textrm{c}_b \textrm{c}_c - \textrm{s}_a \textrm{s}_b) - \textrm{s}_a(\textrm{c}_b\textrm{s}_c + \textrm{c}_c \textrm{s}_b))))^2}{\kappa^2 (\textrm{s}_a^2 (1-(\textrm{c}_a (\textrm{c}_b \textrm{c}_c - \textrm{s}_a \textrm{s}_b) - \textrm{s}_a (\textrm{c}_b\textrm{s}_c + \textrm{c}_c \textrm{s}_b))))}} 
= {\frac{2}{\kappa}\sin s}   
\end{split}
\end{equation*}\\
\noindent In the same way it can be shown that $\cos d(I,B_{[I]}) = \cos d(I,C_{[I]}) = {{(2 \sin s)/\kappa}}.$ Thus, the point $I$ is equidistant to the sides of the triangle and $\cos \textrm{r} = {{(2 \sin s)/\kappa}}.\;\;\Box$\vspace*{1mm}\\

\noindent The equation of the incircle: \noindent When $p_1 := \dfrac 1 {\sin(s{-}a)} , p_2 := \dfrac 1 {\sin(s{-}b)} , p_3 := \dfrac 1 {\sin(s{-}c)}$, 
the equation of the incircle is
\begin{equation*} 
\frac{x_1^2}{p_1^2} + \frac{x_2^2}{p_2^2} + \frac{x_3^2}{p_3^2} - \frac{2x_2 x_3 }{p_2 p_3} - \frac {2x_3 x_1}{p_3 p_1} - \frac{2 x_1 x_2}{p_1 p_2} = 0
\end{equation*}
\noindent \textit{Proof:} This equation is correct because it describes a conic which touches the triangle sides at the traces of $Ge = [p_1:p_2:p_3].\;\;\:\Box$\\\vspace*{-2mm} 

\noindent \textit{Remark}: The incenter of $\Delta_0$ is not only the circumcenter of the tangent triangle but also the circumcenter of the \textit{dual triangle} $\Delta'_{\,0}$. The radius of the circumcircle of the dual triangle $\Delta'_{\,0}$ is $d(I,A') = \arccos \,|\textrm{S}|/\kappa.$

\subsubsection{} \textbf{Triangle centers of $\Delta_1, \Delta_2, \Delta_3$} \hspace*{\fill} \\
\noindent \hspace*{5mm}
The incenter $I$ of $\Delta_0$ can be written $[f(a,b,c):f(b,c,a):f(c,a,b)]$ with $f(a,b,c) = \sin a$. The orthocenter $H$ can also be written
$[f(a,b,c):f(b,c,a):f(c,a,b)]$, but with a different center function $f$; a suitable center function for $H$ is $f(a,b,c) = 1/(\cos b \cos c-\cos a)$.
The first component $f(a,b,c)$ obviously determines the triangle center, the other two one gets by cyclic permutation. \\
Knowing $f(a,b,c)$, we can also write down the corresponding triangle centers for the triangles $\Delta_i, i = 1,2,3:$\\
If $[f(a,b,c):f(b,c,a):f(c,a,b)]$ is representation of a triangle center $Z = Z_0$ of $\Delta_0$ by a barycentric coordinates, then 
\begin{equation*} 
\begin{split}
 Z_1 &= [-f(a_1,b_1,c_1):f(b_1,c_1,a_1):f(c_1,a_1,b_1)] \\
     &= [-f(a,\pi{-}b,\pi{-}c):f(\pi{-}b,\pi{-}c,a):f(\pi{-}c,a,\pi{-}b)]
\end{split}
\end{equation*}
is the corresponding triangle center of $\Delta_1$ and 
\begin{equation*} 
\begin{split}
 Z_2 &= [f(a_2,b_2,c_2):-f(b_2,c_2,a_2):f(c_2,a_2,b_2)],\\
 Z_3 &= [f(a_3,b_3,c_3):f(b_3,c_3,a_3):-f(c_3,a_3,b_3) ]
\end{split}
\end{equation*} 
are the triangle centers of $\Delta_2$ and $\Delta_3$, respectively. \\
\noindent \hspace*{5mm} In the last subsection we presented the radii of the circumcircle and the incircle of $\Delta_0$ as functions of the side lengths and the semiperimeter:
$\textrm{R} = \textrm{R}(a,b,c,s), \textrm{r} = \textrm{r}(a,b,c,s) $. The radii $\textrm{R}_i$ and $\textrm{r}_i$ of the corresponding circles of $\Delta_i$ are: $\textrm{R}_i = \textrm{R}(a_i,b_i,c_i,s_i) , \textrm{r}_i = \textrm{r}(a_i,b_i,c_i,s_i),$ with $s_i = (a_i+b_i+c_i)/2.$\vspace*{1mm}\\
\textit{Remark:} 
The triangle centers $G_i$ and $I_i$, $i = 1,2,3$, are harmonic associates of $G$ and $I$, respectively. 
The orthocenter $H$ is an absolute triangle center: $H = H_1 = H_2 = H_3$.\\

\subsubsection{} \textbf{The staudtian and the area of a triangle} \hspace*{\fill} \\
\noindent\hspace*{5mm}The \textit{staudtian}\footnote{The name is taken  in honor of the geometer von Staudt [1798-1867], who was the first to use this function in spherical geometry, see [4,\,10].} is a function $n$ which assigns each triple of points a real number; the staudtian of $\Delta$ is 
	\[n(\Delta) = \frac 1 2 |\textrm{S}(\Delta)| =  \frac 1 2 |\,det \begin{pmatrix} A^\circ \\ B^\circ\\ C^\circ\end{pmatrix}|. 
	\]
The staudtian has some characteristics of an area. There are equations such as:
	\[ n(\Delta) = \sin a \,\sin b \,\sin \gamma =  \frac 1 2 \,\sin a \,\sin h_a, \, \text{with} \;h_a = d(A, A_H), 
	\]
and for a point $P \in \Delta_0$ the equation: 
	\[ P = [n(BPC): n(CPA) : n(APB)].
\]
 But $n$ lacks the property of additivity. For $P \in \Delta_0$ the inequality holds:
	\[ n(BPC)+ n(CPA) + n(APB) > n(ABC),
\]
and the value of $n(BPC)+ n(CPA) + n(APB)$ takes its maximum for the incenter $I$. (For a proof of the equivalent statement in the hyperbolic plane, see Horv\'ath [10].) \\

The proper triangle area is given by the excess $2\epsilon$; for the triangle $\Delta_0$ this is 
\[2\epsilon(\Delta_0) = \alpha +\beta + \gamma - \pi.
\]
Adding up all the areas of the triangles $\Delta_i, i = 0,1,2,3$, we get $2\pi$ for the area of the whole elliptic plane.

\subsubsection{} \textbf{Isoconjugation} \hspace*{\fill} \\
\noindent\hspace*{5mm}Let $P = [p_1:p_2:p_3]$ be a point not on a sideline of $\Delta_0$ and let $Q = [q_1:q_2:q_3]$ be a point different from the vertices of $\Delta_0$, then the point $R =[p_1q_2q_3:p_2q_3q_1:p_3q_1q_2]$ is called the \textit{isoconjugate} of $Q$ with respect to the pole $P$ or shorter the $P$-isoconjugate of $Q$.\\ Obviously, if $R$ is the {$P$-isoconjugate} of $Q$, then $Q$ is the {$P$-isoconjugate} of $R$. \vspace*{-3mm}\\

\noindent\textit{Some examples:}\\
\hspace*{4mm}$G$-isoconjugation is also called \textit{isotomic conjugation}. The fixed points of $G$-isoconjugation are the centroids $G, G_1, G_2, G_3$.\\
The isotomic conjugate of the Gergonne point $Ge$ is the \textit{Nagel point} $N\!a$. Its traces are the touch points of the excircles of $\Delta_0$ (= incircles of the triangles $\Delta_i, i = 1,2,3$) with the sides of $\Delta_0$.

Isogonal conjugation leaves the incenters $I_i, i = 0,1,2,3$, fixed. It is the $P$-isoconjugation for $P = [\sin^2 a:\sin^2 b:\sin^2 c]$. This point is called \textit{symmedian} and is usually notated by the letter $K$.\\
The circumcenters $O_1,O_2,O_3$ form a perspective triple. The perspector is the isogonal conjugate of $O$, and we denote this point by $H^-$.

Define the point $\tilde{K}$ by\\
\centerline{$\tilde{K} := [1{-}\cos a, 1{-}\cos b, 1{-}\cos c,] = [\sin^2 a/2, \sin^2 b/2,\sin^2 c/2].$}\\
Horv\'ath uses the name \textit{Lemoine point} for it and we also shall use this name. \\
The points on the circumcircle: {$(1-\cos a) x_2 x_3 + (1-\cos b) x_3 x_1 + (1-\cos c) x_1 x_2 = 0$} are the $\tilde{K}$-isoconjugates of the points on the tripolar of $G$:
${x_1 + x_2 +  x_3 = 0},$ a line which is also the dual of $O$.\\

\subsection{} \textbf{Conics} \hspace*{\fill}\vspace*{-1mm} 
\subsubsection{} \textbf{Different types of conics} \hspace*{\fill} \\
\noindent\hspace*{5mm}Let $\mathfrak{M} = (\mathfrak{m}_{ij})_{i,j=1,2,3}$ be a symmetric matrix, then 
the quadratic equation
	\[ \mathfrak{m}_{11}x_1^2 + \mathfrak{m}_{22}x_2^2 + \mathfrak{m}_{33}x_3^2 + 2\mathfrak{m}_{23}x_2 x_3 + 2\mathfrak{m}_{31}x_3 x_1 + 2\mathfrak{m}_{12}x_1 x_2 = 0
\]
is the equation of a conic which we denote by $\mathcal{C}(\mathfrak{M})$. Given a symmetric matrix $\mathfrak{M}$ and a nonzero real number $t$, then the conics $\mathcal{C}(\mathfrak{M})$ and $\mathcal{C}(t\mathfrak{M})$ are the same. We can reverse this implication if we restrict ourselves to real conics. Here, we define: A \textit{real conic} in $\mathcal{P}$ is
\begin{itemize} 
\item[-] either the union of two different real lines,
\item[-] or a circle (with radius $r \in [0,\frac{\pi}{2}]$) ,
\item[-] or a proper ellipse, that is an irreducible ($det(\mathfrak{M}) \ne 0$) conic with infinitely many real points and which is not a circle.
\end{itemize} 
\textit{Remarks}: 
\begin{itemize} 
\item[-]A double line can be regarded as a circle with radius $r = \frac{\pi}{2}$.
\item[-]There is no difference between ellipses, hyperbolae and parabolae in elliptic geometry, cf. [8].
\end{itemize} 

The polar of a point $P = [p_1:p_2:p_3] $ with respect to the conic $\mathcal{C}(\mathfrak{M})$ is the line with the equation\vspace*{-2mm}
	\[ (x_1,x_2,x_3)\,\mathfrak{M}\!\begin{pmatrix} p_1 \\ p_2\\ p_3 \end{pmatrix} = 0.
\]
The pole of the line $\ell: \ell_1 x_1 + \ell_2 x_2 + \ell_3 x_3 = 0$ is the point $P = [p_1:p_2:p_3]$ with
	$(p_1,p_2,p_3) = (\ell_1,\ell_2, \ell_3)\mathfrak{M}^{\#}$\vspace*{1mm}, 
where\\ 
\centerline{$\mathfrak{M}^\# = \begin{pmatrix} \mathfrak{m}_{22}\mathfrak{m}_{33}-\mathfrak{m}_{23}^2&\mathfrak{m}_{23}\mathfrak{m}_{31}-\mathfrak{m}_{12}\mathfrak{m}_{33}&{\mathfrak{m}_{12}\mathfrak{m}_{23}-\mathfrak{m}_{31}\mathfrak{m}_{22}}\\ 
\mathfrak{m}_{23}\mathfrak{m}_{31}-\mathfrak{m}_{12}\mathfrak{m}_{33}&\mathfrak{m}_{33}\mathfrak{m}_{11}-\mathfrak{m}_{31}^2&\mathfrak{m}_{31}\mathfrak{m}_{12}-\mathfrak{m}_{13}\mathfrak{m}_{11} \\
\mathfrak{m}_{12} \mathfrak{m}_{23}-\mathfrak{m}_{31}\mathfrak{m}_{22}&\mathfrak{m}_{31}\mathfrak{m}_{12}-\mathfrak{m}_{13}\mathfrak{m}_{11}&\mathfrak{m}_{11}\mathfrak{m}_{22}-\mathfrak{m}_{12}^2 \end{pmatrix}$}\vspace*{1.5mm}\\  
is the adjoint of $\mathfrak{M}.$\\

\subsubsection{} \textbf{The perspector of a conic} \hspace*{\fill} \\
\noindent\hspace*{5mm}If $\mathfrak{M}$ is a diagonal matrix, then the polar lines of $A, B, C$ are the lines $B\vee C$, $C\vee A$, $A\vee B$, respectively. If $\mathfrak{M}$ is not diagonal,
then the poles of $B\vee C$, $C\vee A$, $A\vee B$ with respect to the conic form a perspective triple with perspector 
	\[[\frac 1 {\mathfrak{m}_{11} \mathfrak{m}_{23}- \mathfrak{m}_{31} \mathfrak{m}_{12}}: \frac 1 {\mathfrak{m}_{22} \mathfrak{m}_{31}- \mathfrak{m}_{12} \mathfrak{m}_{23}}:\frac 1 {\mathfrak{m}_{33} \mathfrak{m}_{12}- \mathfrak{m}_{23} \mathfrak{m}_{31}}].
\]
\noindent\textit{An example}: A matrix of the absolute conic is $\mathfrak{T}$ and the perspector is $H$. Since there are no real points on this conic, it is not a real conic.\\

\subsubsection{} \textbf{Symmetry points and symmetry axes of real conics} \hspace*{\fill} \\
\noindent\hspace*{5mm}A point $P = [p_1:p_2:p_3]$ is a \textit{symmetry point} of a conic if for every point $Q = [q_1:q_2:q_3]$ on this conic the mirror image of $Q$ in $P$ is also a point on this conic. A line $l$ is a \textit{symmetry axis} if its dual $l^\delta$ is a symmetry point. The meet of two different symmetry axes is a symmetry point and the join of two different symmetry points is a symmetry axis.\vspace*{1mm}\vspace*{2mm}\\
\noindent\textit{Three examples}:\\ 
\hspace*{5mm}If the conic is the union of two different lines which meet at a point $P$, then the point $P$ and the duals of the two angle bisectors are the symmetry points. The point $P$ is regarded as the center of the conic. \\
\hspace*{5mm}The symmetry points of a circle with center $P$ and radius $r \in [0,\pi/2 ]$ are the center $P$ and all points on $P^\delta$.\\ 
\hspace*{5mm}A proper ellipse has three symmetry points, one lies inside the ellipse and is regarded as its center, the other two lie outside.\\

How can we find the symmetry points of a real conic $\mathcal{C}(\mathfrak{M})$ in case of a circle with radius $r \in\, ]0,\pi/2 [$ or a proper ellipse? In this case, a point $P = {[p_1:p_2:p_3]}$ is a symmetry point precisely when its dual $P^\delta$ is identical with the polar of $P$ with respect to the conic, this is when the vector $(p_1,p_2,p_3)$ is an eigenvector of the matrix $\mathfrak{T}^{\#}\mathfrak{M}$.\vspace*{2mm}\\

\noindent\textit{Three examples}:\\ 
\hspace*{5mm}A circumconic with perspector $P = {[p_1:p_2:p_3]}$ is described by the equation\\
\centerline{$p_1 x_2 x_3 + p_2 x_3 x_1 + p_3 x_1 x_2 = 0$.}\\
A comparison of this equation with the equation of the circumcircle shows that the Lemoine point $\tilde{K}$ is the perspector of the circumcircle. The symmetry points are, besides the circumcenter $O$, the points on the line $G^\tau: x_1 + x_2 + x_3 = 0.$\\
\hspace*{5mm}We want to determine the symmetry points of the circumconic in case of $P = {[1+ 2\textrm{c}_a:1+ 2\textrm{c}_b:1+ 2\textrm{c}_c]}.$ If the triangle $\Delta_0$ is not equilateral, then the circumconic is a proper ellipse  and the matrix $\mathfrak{T}^{\#}\mathfrak{M}$ has three different eigenvalues. One is $2 (\textrm{c}_a^2  + \textrm{c}_b^2 + \textrm{c}_c^2 - 2 \textrm{c}_a \textrm{c}_b \textrm{c}_c - 1)$, belonging to the eigenvector $(1,1,1)$. The other two eigenvalues and their corresponding eigenvectors can also be explicitly calculated. Here, the characteristic polynomial of  $\mathfrak{T}^{\#}\mathfrak{M}$ splits into a linear and a quadratic rational factor. But for proper ellipses this is an exception. In general, formulas for the symmetry points of a circumconic with a given perspector are rather complicated.\\ On the other hand, knowing the center $M = [m_1: m_2:m_3]$ of a circumconic, its perspector $P$ can be calculated quite easily:\vspace*{-2mm} 
	\[
\begin{split}	
	P = [m_1(2m_2 m_3 \textrm{c}_a - m_1^2+m_2^2+m_3^2)&: m_2(2m_3 m_1 \textrm{c}_b + m_1^2-m_2^2+m_3^2)\\
	&: m_3(2m_1 m_2 \textrm{c}_c + m_1^2+m_2^2-m_3^2)].
\end{split}	
\]
\hspace*{5mm}The equation of the inconic with perspector $P = {[p_1:p_2:p_3]}$ is
	\[\frac {x_1^2}{p_1^2}+\frac {x_2^2}{p_2^2}+\frac {x_3^2}{p_3^2}+2 \frac {x_2 x_3} {p_2 p_3}+2 \frac {x_3 x_1} {p_3 p_1}+2 \frac {x_1 x_2} {p_1 p_2} = 0
\]
\noindent If the perspector $P$ is the Gergonne point $Ge$, then the inconic is the incircle. Its symmetry points are the incenter $I$ and the points on the line
	\[ 
\begin{split}		
	\textrm{s}_{s-a}(\textrm{s}_{a}\textrm{s}_{b}\textrm{s}_{c} - \textrm{s}_{a} \textrm{S}_{A} + \textrm{s}_{b} \textrm{S}_{B} + \textrm{s}_{c} \textrm{S}_{C}) x_1 &+ \textrm{s}_{s-b}(\textrm{s}_{a}\textrm{s}_{b}\textrm{s}_{c} + \textrm{s}_{a} \textrm{S}_{A} - \textrm{s}_{b} \textrm{S}_{B} + \textrm{s}_{c} \textrm{S}_{C}) x_2 \\
	&+ \textrm{s}_{s-c}(\textrm{s}_{a}\textrm{s}_{b}\textrm{s}_{c} + \textrm{s}_{a} \textrm{S}_{A} + \textrm{s}_{b} \textrm{S}_{B} - \textrm{s}_{c} \textrm{S}_{C}) x_3 = 0.
\end{split}			
\]

\subsubsection{} \textbf{Bicevian conics} \hspace*{\fill} \\
\noindent\hspace*{5mm}Let $P = [p_1:p_2:p_3]$ and $Q = [q_1:q_2:q_3]$ be two different points, not on any of the lines of $\Delta_0$. Then \vspace*{1 mm}
	\[ (x_1,x_2,x_3) \begin{pmatrix}\frac 2{p_1 q_1}& \frac 1{p_1 q_2} + \frac 1{p_2 q_1}&\frac 1{p_3 q_1} + \frac 1{p_1 q_3}\\  \frac 1{p_1 q_2} +\frac 1{p_2 q_1}& \frac 2{p_2 q_2}& \frac 1{p_2 q_3} + \frac 1{p_3 q_2} \\ \frac 1{p_3 q_1} + \frac 1{p_1 q_3}& \frac 1{p_2 q_3} + \frac 1{p_3 q_2}  & \frac 2{p_3 q_3} \end{pmatrix}\!\begin{pmatrix} x_1 \\ x_2\\ x_3 \end{pmatrix} = 0\vspace*{1 mm}
\]
is the equation of a conic, which passes through the traces of $P$ and $Q$. The perspector of this conic is the point $R = [r_1:r_2:r_3]$ with 
	\[ (r_1,r_2,r_3) = (p_1 q_1 (p_2 q_3 - p_3 q_2), p_2 q_2 (p_3 q_1 - p_1 q_3), p_3 q_3 (p_1 q_2 - p_2 q_1)).
\]

What are the conditions for $Q$ to be a circumcevian conjugate of $P$ or, in other words, what are the conditions for the bicevian conic to be a circle? Here, we should keep in mind that the point triple $A_PB_PC_P$ has four circumcircles and, therefore, there are four circumcevian conjugates of $P$. First, we calculate the barycentric coordinates of the circumcenter $M_0$ of the triangle $(A_P B_P C_P)_0$: The dual $M_0^{\;\delta}$ of $M_0$ is the line $(A_P - B_P) \vee (A_P - C_P)$, which is described by the equation
	\[  
	p_2 p_3(-t_1+t_2 +t_3) x_1 + p_2 p_1(t_1-t_2 +t_3) x_2 + p_1 p_2(t_1+t_2 -t_3) x_3 = 0, 
	\]																											
with\; $t_1:= \|(0,p_2,p_3)\|_\ast = \sqrt{p_2^{\;2}+2 p_2 p_3 \textrm{c}_a+p_3^{\;2}}$,\\
\hspace*{7.5mm}	$t_2:= \|(p_1,0,p_3)\|_\ast = \sqrt{p_3^{\;2}+2 p_3 p_1 \textrm{c}_b+p_1^{\;2}}$ \\
and $\;\, t_3:= \|(p_1,p_2,0)\|_\ast =\sqrt{p_1^{\;2} + 2 p_1 p_2 \textrm{c}_c + p_2^{\;2}}.$\\
We put $(s_1,s_2,s_3) :=  ((-t_1+t_2 +t_3)/p_1,(t_1-t_2 +t_3)/p_2,(t_1+t_2 -t_3)/p_3)$
and calculate the coordinates $[m_{01}:m_{02}:m_{03}]$ of $M_0$ and the coordinates $[q_{01}:q_{02}:q_{03}]$ of the  circumcevian conjugate $Q_0$ of $P$:
	\[
			\begin{split}
		(m_{01},m_{02},m_{03}) &= (s_1,s_2,s_3) \mathfrak{T}^{\#}\\
		                       &= ((-t_1+t_2 +t_3)/p_1,(t_1-t_2 +t_3)/p_2,(t_1+t_2 -t_3)/p_3) \mathfrak{T}^{\#}\\
		(q_{01},q_{02},q_{03})\;\;\; &= \big(\frac {1}{(s_1^{\,2} - 4)p_1}, \frac {1}{(s_2^{\,2} - 4)p_2},\frac {1}{(s_3^{\,2} - 4)p_3}\big)\\
		&= \big(\frac {p_1}{({-}t_1{+}t_2{+}t_3)^2-4 p_1^{\,2}}, \frac {p_2}{(t_1{-}t_2{+}t_3)^2-4 p_2^{\,2}},\frac {p_3}{(t_1{-}t_2{+}t_3)^2-4 p_3^{\,2})}\big) .
\end{split}	
\]
For the radius $r$ of the cevian-circle we calculate 
	\[r = \arccos \frac {2 \textrm{S}^2}{\|(m_{01},m_{02},m_{03})\|_\ast}.
\]
\noindent \hspace*{5mm}The other circumcenters $M_i$ and circumcevian conjugates $Q_i$ , $i = 1,2,3$, can be determined in the same way. We give the results for $i = 1$:
\[
		\begin{split}
		&(m_{11},m_{12},m_{13}) = ((t_1{+}t_2{+}t_3)/p_1,({-}t_1{-}t_2{+}t_3)/p_2,({-}t_1{+}t_2{-}t_3)/p_3) \mathfrak{T}^{\#},\\
		&(q_{11},q_{12},q_{13}) = \big(\frac {p_1}{(t_1{+}t_2{+}t_3)^2-4 p_1^{\,2}}, \frac {p_2}{({-}t_1{-}t_2{+}t_3)^2-4 p_2^{\,2}},\frac {p_3}{({-}t_1{+}t_2{-}t_3)^2-4 p_3^{\,2})}\big).
	\end{split}		
\]

\section{Four central lines}
\subsection{} \textbf{The orthoaxis $G^+\! \vee H$} \hspace*{\fill} 
\subsubsection{} \textbf{Triangle centers on the orthoaxis} \hspace*{\fill} \\
\noindent\hspace*{5mm}Wildberger [27] introduces the name \textit{orthoaxis} for a line incident with several triangle centers: the orthocenter $H$, the centroid $G^+$ of the antimedial triangle, the orthostar $H^{\star}$ and three more triangle centers:
\begin{itemize} 
\item[-]One is the point $O^+ := [\textrm{S}_A \textrm{s}_a^{\,2}:\textrm{S}_B \textrm{s}_b^{\,2}:\textrm{S}_C \textrm{s}_c^{\,2}]$, the isogonal conjugate of $H$.\\
In [27, 28] it is called \textit{basecenter} and introduced as the meet of the lines $A\vee (B_H \vee C_H)^\delta$, $B\vee (C_H \vee A_H)^\delta,$  $C\vee (A_H \vee B_H)^\delta.$ In [26] the point $O^+$ is called \textit{pseudo-circumcenter}, and it is shown that it is the meet of the perpendicular pseudo-bisectors $A_{G^+}\!  \vee A'$, $B_{G^+}\! \vee B',$ $C_{G^+}\! \vee C'$.
\item[-]The second point is \\
$L :=  [{\textrm{c}_a·(-\textrm{c}_a^{\,2}{+}\textrm{c}_b^{\,2}{+}\textrm{c}_c^{\,2}{+}1){-}2\textrm{c}_b\textrm{c}_c}:{\textrm{c}_b(\textrm{c}_a^{\,2}{-}\textrm{c}_b^{\,2}{+}\textrm{c}_c^{\,2}{+}1){-}2\textrm{c}_a\textrm{c}_c}\\
\hspace*{51mm}: {\textrm{c}_c(\textrm{c}_a^{\,2}{+}\textrm{c}_b^{\,2}{-}\textrm{c}_c^{\,2}{+}1){-}2\textrm{c}_a\textrm{c}_b}].$\\
It is the common point of the lines $A^{G^+}\!\! \vee A'$, $B^{G^+}\!\! \vee B'$, $C^{G^+}\!\! \vee C'$, and it is called \textit{double dual point} in [27, 28].
\item[-] The third point is the intersection of the orthoaxis with the orthic axis $H^\tau$ and has coordinates\\
	$[\textrm{S}_B \textrm{S}_C (-2 \textrm{c}_a^2 + \textrm{c}_b^2 + \textrm{c}_c^2): \textrm{S}_C \textrm{S}_A (\textrm{c}_a^2 - 2 \textrm{c}_b^2 + \textrm{c}_c^2):\textrm{S}_A \textrm{S}_B (\textrm{c}_a^2 + \textrm{c}_b^2 - 2 \textrm{c}_c^2)].$\vspace*{2mm}
\end{itemize}
\subsubsection{} \textbf{The bicevian conic through the traces of $H$ and $G^+$} \hspace*{\fill} \\
\noindent\hspace*{5mm}We prove a conjecture of Vigara [26]:\\  
The orthoaxis $\ell = H  \vee G^+$ is a symmetry axis of the bicevian conic which passes through the traces of $H$ and $G^+$. Besides $\ell^\delta$, the points $O^+ \!+ H$, $O^+ \!- H$ are symmetry points of this conic.\footnote{The statement in the previous version that these last two points are in general not symmetry points is wrong!}   \\ \vspace*{-2mm}

\noindent \textit{Proof:}
The orthoaxis is described by the equation: 
	\[
	\textrm{S}_A (\textrm{c}_b^{\,2} - \textrm{c}_c^{\,2})x_1 
	+ \textrm{S}_B (\textrm{c}_c^{\,2} - \textrm{c}_a^{\,2})x_2 
  +\textrm{S}_C (\textrm{c}_a^{\,2} - \textrm{c}_b^{\,2})x_3 = 0,
\]
so its dual is the point $P:= \ell^\delta = [\textrm{c}_a (\textrm{c}_b^{\,2} - \textrm{c}_c^{\,2}), \textrm{c}_b (\textrm{c}_c^{\,2} - \textrm{c}_a^{\,2}), \textrm{c}_c  (\textrm{c}_a^{\,2} - \textrm{c}_b^{\,2})]$.\\
This point $P$ is also the perspector of the conic through the traces of $H$ and $G^+$; the equation of the conic is \vspace*{-2mm}
\[ (x_1,x_2,x_3)\, {\mathfrak{M}}\!\begin{pmatrix} x_1 \\ x_2\\ x_3 \end{pmatrix} = 0,\vspace*{-3mm}
\]
with 
	\[\mathfrak{M} = \begin{pmatrix}2\,\textrm{S}_A \textrm{c}_b \textrm{c}_c &-\textrm{c}_c(\textrm{S}_A \textrm{c}_a + \textrm{S}_B \textrm{c}_b) &-\textrm{c}_b(\textrm{S}_C \textrm{c}_c + \textrm{S}_A \textrm{c}_a)\\ -\textrm{c}_c(\textrm{S}_A \textrm{c}_a + \textrm{S}_B \textrm{c}_b)  & 2\,\textrm{S}_B \textrm{c}_c \textrm{c}_a &-\textrm{c}_a(\textrm{S}_B \textrm{c}_b + \textrm{S}_C \textrm{c}_c) \\ -\textrm{c}_b(\textrm{S}_C \textrm{c}_c + \textrm{S}_A \textrm{c}_a)& -\textrm{c}_a(\textrm{S}_B \textrm{c}_b + \textrm{S}_C \textrm{c}_c) & 2\,\textrm{S}_C \textrm{c}_a \textrm{c}_b \end{pmatrix}.
\]
The point $P$ is not only the perspector of the conic but also a symmetry point: The polar of $P$ with respect to this conic is calculated by 
\[ (x_1,x_2,x_3)\, {\mathfrak{M}}\!\begin{pmatrix} \textrm{c}_a (\textrm{c}_b^{\,2} - \textrm{c}_c^{\,2}) \\  \textrm{c}_b (\textrm{c}_c^{\,2} - \textrm{c}_a^{\,2})\\\textrm{c}_c  (\textrm{c}_a^{\,2} - \textrm{c}_b^{\,2}) \end{pmatrix} = 0.
\]
and this, again, is an equation of the orthoaxis. Since  $P$ is a symmetry point, the orthoaxis is a symmetry line of the conic.\\
Now we will show that the points $O^+ \!+ H$, $O^+ \!- H$ on $\ell$ are the other two symmetry points of $\mathcal{C}(\mathfrak{M})$. It is sufficient to show the correctness of the equation $(H^\circ+O^{+\circ})\mathfrak{M}\cdot(H^\circ-O^{+\circ})=0$, which expresses that $H\mp O^+$ is a point on the polar of $H\pm O^+$ with respect to $\mathcal{C}(\mathfrak{M})$. We transform the equation equivalently:\vspace*{-3mm}\\
 \[
\begin{split}		
&\;\;\;\;\;\;\;\;\;\;\;\;\;\,\,(H^\circ + O^{+\circ})\mathfrak{M}\cdot(H^\circ - O^{+\circ}) =0\\
\Leftrightarrow\;\;\;&\;\;\;\;\;\;\;\;\;\;\;\;\;\;\;\;\;\;H^\circ\mathfrak{M}\cdot H^\circ - O^{+\circ}\mathfrak{M} \cdot O^{+\circ} =0\\
\Leftrightarrow\;\;\;&(\boldsymbol{h}\,\mathfrak{M}\cdot \boldsymbol{h})(\boldsymbol{o}\,\mathfrak{T}\cdot \boldsymbol{o}) - (\boldsymbol{h}\,\mathfrak{T}\cdot \boldsymbol{h})(\boldsymbol{o}\,\mathfrak{M}\cdot \boldsymbol{o}) =0\\ 
&\;\text{with}\; \boldsymbol{h} = (\textrm{S}_B \textrm{S}_C,\textrm{S}_C \textrm{S}_A,\textrm{S}_A \textrm{S}_B)\; \textrm{and}\;\boldsymbol{o} = (\textrm{S}_A \textrm{s}_a^{\,2},\textrm{S}_B \textrm{s}_b^{\,2},\textrm{S}_C \textrm{s}_c^{\,2}).
\vspace*{-2mm}
\end{split}		
\]
The last equation holds as can be checked with the help of a CAS.$\;\Box$

\subsection{} \textbf{The line $G \vee O$} \hspace*{\fill} 

\subsubsection{} \textbf{Triangle centers on the line $G \vee O$} \hspace*{\fill} \\
\noindent\hspace*{5mm}The line $G \vee O$ is the orthoaxis of the medial triangle and has the equation 
	\[
\begin{split}		
	(1+ \textrm{c}_a - \textrm{c}_b - \textrm{c}_c) (\textrm{c}_b - \textrm{c}_c)x_1 &+  (1 -	\textrm{c}_a + \textrm{c}_b - \textrm{c}_c ) (\textrm{c}_c - \textrm{c}_a)x_2\\
	&+ (1- \textrm{c}_a -\textrm{c}_b + \textrm{c}_c)(\textrm{c}_a - \textrm{c}_b)x_3 = 0.
\end{split}		
\]
Besides $G$ and $O$ it contains the following triangle centers: 
\begin{itemize} 
\item[-] The isogonal conjugate of $O$ with barycentric coordinates 
\[[\frac{1+\textrm{c}_a}{1+ \textrm{c}_a - \textrm{c}_b - \textrm{c}_c}:\frac{1+\textrm{c}_b}{1- \textrm{c}_a + \textrm{c}_b - \textrm{c}_c}:\frac{1+\textrm{c}_c}{1- \textrm{c}_a - \textrm{c}_b + \textrm{c}_c}].
\]
This point is also the common point of the lines ${\textrm{perp}(B_G \vee C_G,A)}$,\\ ${\textrm{perp}(C_G \vee A_G,B)}$ and
${\textrm{perp}(A_G \vee B_G,C)}$; we denote it by $H^-$.\vspace*{-1.5mm}\\
\item[-] The $\tilde{K}$-conjugate of $O$; its coordinates are 
\[[\frac{1}{1+ \textrm{c}_a - \textrm{c}_b - \textrm{c}_c}:\frac{1}{1- \textrm{c}_a + \textrm{c}_b - \textrm{c}_c}:\frac{1}{1- \textrm{c}_a - \textrm{c}_b + \textrm{c}_c}].
\]
\item[-] The point $L$  was already introduced in the last subsection; it is the intersection of the line $G\vee O$ with the orthoaxis. But it is also  a point on the line $I \vee Ge$ and a point on all the lines $ G_i \vee O_i, i = 0,1,2,3$. We call this point \textit{de\,Longchamps point}. The main reason for choosing the name is: This point $L$ is the radical center of the three power circles $\mathcal{C}(B_G + C_G, A), \mathcal{C}(C_G + A_G, B), \mathcal{C}(A_G + B_G, C).$\vspace*{1mm}
\end{itemize}
\textit{Proof} of the last statement: We will outline the proof for $\textrm{perp}(B_G\vee C_G,A)$ being the radical line of the first two power circles. In a similar way it can be shown that $\textrm{perp}(C_G \vee A_G,B), \textrm{perp}(A_G \vee B_G,C)$ are the other two radical lines.\\ We start with the dual of the line through the centers of the first two circles; this is the point
	\[P = [(1{-}\textrm{c}_a) (1{+}\textrm{c}_a{+}\textrm{c}_b{+}\textrm{c}_c): -(1{+}\textrm{c}_b)(1{-}\textrm{c}_a{-}\textrm{c}_b{+}\textrm{c}_c): -(1{+}\textrm{c}_c)(1{-}\textrm{c}_a{+}\textrm{c}_b{-}\textrm{c}_c)].
\]
Then for every real number $t$, the vector 
	\[
\begin{split}	
	\boldsymbol{p}_t := (&(1 - \textrm{c}_a) (1+ \textrm{c}_a + \textrm{c}_b + \textrm{c}_c) +t({\textrm{c}_a·(1{-}\textrm{c}_a^{\,2}{+}\textrm{c}_b^{\,2}{+}\textrm{c}_c^{\,2}){-}2\textrm{c}_b\textrm{c}_c}) ,\\ -&(1+ \textrm{c}_b)(1 - \textrm{c}_a - \textrm{c}_b + \textrm{c}_c)+t({\textrm{c}_b(1{+}\textrm{c}_a^{\,2}{-}\textrm{c}_b^{\,2}{+}\textrm{c}_c^{\,2}){-}2\textrm{c}_a\textrm{c}_c}),\\ -&(1+\textrm{c}_c)(1-\textrm{c}_a +\textrm{c}_b - \textrm{c}_c)+t({\textrm{c}_c(1{+}\textrm{c}_a^{\,2}{+}\textrm{c}_b^{\,2}{-}\textrm{c}_c^{\,2}){-}2\textrm{c}_a\textrm{c}_b})).
\end{split}
\]
represents the barycentric coordinates of a point on $P\vee L.$ We substitute the components of $\boldsymbol{p}_t$ for $x_1,x_2,x_3$ in the equations of the two power circles
	\[
\begin{split}	
(x_1 (\textrm{c}_b{+}\textrm{c}_c){+}(\textrm{c}_a{+}1) (x_2{+}x_3))^2&=(2 \textrm{c}_a x_2 x_3{+}2 \textrm{c}_b x_1 x_3{+}2 \textrm{c}_c x_1 x_2{+}x_1^2{+}x_2^2{+}x_3^2)(\textrm{c}_b{+}\textrm{c}_c)^2,\\
(x_2 (\textrm{c}_a{+}\textrm{c}_c){+}(\textrm{c}_b{+}1) (x_1{+}x_3))^2&=(2 \textrm{c}_a x_2 x_3{+}2 \textrm{c}_b x_1 x_3{+}2 \textrm{c}_c x_1 x_2{+}x_1^2{+}x_2^2{+}x_3^2) (\textrm{c}_a{+}\textrm{c}_c)^2,	
\end{split}	
\]
 solve for $t$ and get the same solutions for both circles.$\;\;\Box$ \\

\noindent The points $O, G, H^-, L$ form a harmonic range.\\ \vspace*{-2 mm}

\noindent \textit{Proof:} We introduce the vectors \\
$\boldsymbol{o} = ((1{-}\textrm{c}_a)(1{+}\textrm{c}_a{-}\textrm{c}_b{-}\textrm{c}_c) , (1{-}\textrm{c}_b)(1{-}\textrm{c}_a{+}\textrm{c}_b{-}\textrm{c}_c), (1{-}\textrm{c}_c)(1{-}\textrm{c}_a{-}\textrm{c}_b{+}\textrm{c}_c))$, \\
$\boldsymbol{g} = (1,1,1)$,  \\
$\boldsymbol{h} =((1{+}\textrm{c}_a)/(1{+}\textrm{c}_a{-}\textrm{c}_b{-}\textrm{c}_c) , (1{+}\textrm{c}_b)/(1{-}\textrm{c}_a{+}\textrm{c}_b{-}\textrm{c}_c), (1{+}\textrm{c}_c)/(1{-}\textrm{c}_a{-}\textrm{c}_b{+}\textrm{c}_c))$,\\
$\boldsymbol{l}=({\textrm{c}_a·(-\textrm{c}_a^{\,2}{+}\textrm{c}_b^{\,2}{+}\textrm{c}_c^{\,2}{+}1){-}2\textrm{c}_b\textrm{c}_c},{\textrm{c}_b(\textrm{c}_a^{\,2}{-}\textrm{c}_b^{\,2}{+}\textrm{c}_c^{\,2}{+}1){-}2\textrm{c}_a\textrm{c}_c}, {\textrm{c}_c(\textrm{c}_a^{\,2}{+}\textrm{c}_b^{\,2}{-}\textrm{c}_c^{\,2}{+}1){-}2\textrm{c}_a\textrm{c}_b})$,\\
the real numbers\\
$r =  1 - \textrm{c}_a^2 - \textrm{c}_b^2 - \textrm{c}_c^2 + 2 \textrm{c}_a \textrm{c}_b \textrm{c}_c,$ \\
$s =  (1+\textrm{c}_a - \textrm{c}_b - \textrm{c}_c) (1-\textrm{c}_a + \textrm{c}_b - \textrm{c}_c) (1 -\textrm{c}_a - \textrm{c}_b + \textrm{c}_c)$,\\
$t = (1+ \textrm{c}_a + \textrm{c}_b + \textrm{c}_c)$,\\
and get the equations
	$\;\;s \boldsymbol{o} + t\boldsymbol{h} = 2r\boldsymbol{g}\;\;$ and $\;\; s \boldsymbol{o} - t\boldsymbol{h} = 2\boldsymbol{l}\;\;\;\;\;\;\Box.$

\subsubsection{} \textbf{A cubic curve as a substitute for the Euler circle} \hspace*{\fill} \\
\noindent\hspace*{5mm}
For a point $P = [p_1:p_2:p_3] \in \mathcal{P}$  we calculate the pedals of $P$ on the sidelines of the medial triangle:
	\[ \tilde{A}_{[P]} := {(P \vee (B_{[G]} \vee C_{[G]})^\delta) \wedge (B_{[G]} \vee C_{[G]}}), \tilde{B}_{[P]} , \tilde{C}_{[P]}.
\]
The points $P$ for which $\tilde{A}_{[P]}, \tilde{B}_{[P]}, \tilde{C}_{[P]}$ are collinear lie on a cubic that passes through the traces of the points $H^-$ and $G_i, i = 0,1,2,3.$ We call this cubic \textit{Euler-Feuerbach cubic}. In metric affine geometries this cubic splits into the nine-point-circle (Euler circle) and the line $G^\tau\!: x_1+x_2+x_3 = 0$. \\ \vspace*{-2 mm}

\noindent \textit{Proof:} Instead of proving the statement for the triangle $\Delta_0$, we present the proof for the antimedial triangle of $\Delta_0$. In this case, the formulae obtained are substantially shorter. \\
For the coordinates of the pedals $A_{[P]}, B_{[P]}, C_{[P]}$ of a point $P$ on the sidelines of $\Delta$, see 2.3.9. If these three pedals are collinear on a line, then this line is called a \textit{Simson line} of $P$. The locus of points $P$ having a Simson line is the cubic with the equation
	\[
	\begin{split}	
	\textrm{c}_a x_1(x_2^{\,2}\textrm{s}_c^{\,2} + x_3^{\,2}\textrm{s}_b^{\,2}) + \textrm{c}_b x_2(x_3^{\,2}\textrm{s}_a^{\,2} + x_1^{\,2}\textrm{s}_c^{\,2}) &+ \textrm{c}_c x_3(x_1^{\,2}\textrm{s}_b^{\,2} + x_2^{\,2}\textrm{s}_a^{\,2}) \\
	&- 2 x_1 x_2 x_3 (1- \textrm{c}_a \textrm{c}_b \textrm{c}_c) = 0.
\end{split}	
\]
The centroid of the antimedial triangle is $G^+ =[\textrm{c}_a:\textrm{c}_b:\textrm{c}_c]$. The traces of $G^+$  on the sidelines of the antimedial triangle are $A, B, C$. The tripolar of $G^+$ meets the sidelines of the antimedial triangle in $[0: \textrm{c}_b: -\textrm{c}_c], [-\textrm{c}_a:0:\textrm{c}_c]$ and $[\textrm{c}_a:-\textrm{c}_b:0]$.\vspace*{1mm}
The traces of the point $L$ on the sidelines of the antimedial triangle are 
	\[
		\begin{split}	
		&[ \textrm{c}_a (\textrm{c}_b^{\,2}+\textrm{c}_c^{\,2}) - 2 \textrm{c}_b^{\,2} \textrm{c}_b^{\,2}:  \textrm{c}_b (\textrm{c}_b^{\,2} - \textrm{c}_c^{\,2}): \textrm{c}_c (\textrm{c}_b^{\,2} - \textrm{c}_c^{\,2})],\\
	&[ \textrm{c}_a (\textrm{c}_c^{\,2} - \textrm{c}_a^{\,2}): \textrm{c}_b (\textrm{c}_c^{\,2}+\textrm{c}_a^{\,2}) - 2 \textrm{c}_b^{\,2} \textrm{c}_b^{\,2}: \textrm{c}_c (\textrm{c}_c^{\,2} - \textrm{c}_a^{\,2})],	\\
	&[ \textrm{c}_a (\textrm{c}_a^{\,2} - \textrm{c}_b^{\,2}):  \textrm{c}_b (\textrm{c}_a^{\,2} - \textrm{c}_b^{\,2}): \textrm{c}_c (\textrm{c}_a^{\,2}+\textrm{c}_b^{\,2}) - 2 \textrm{c}_a^{\,2} \textrm{c}_b^{\,2}].	
\end{split}		
\]
It can now be checked that the coordinates of all the traces satisfy the cubic equation. $\;\;\Box$\vspace*{1mm}\\
\textit{Remarks:}\\
This circumcubic of $ABC$ is the non pivotal isocubic $n\mathcal{K}(K,G^+,t)$ with pole $K$ (symmedian), root $G^+$ and parameter $t = {-2(1- \textrm{c}_a \textrm{c}_b \textrm{c}_c)}{\ne 0}$ (terminology adopted from [5, 7]).
It is the locus of dual points of the Simson lines of $\Delta$, and it also passes through 
\begin{itemize} 
\item[-] the vertices $A', B', C'$ of the dual triangle,
\item[-] and through the points $[- \textrm{S}_B \textrm{S}_C : \textrm{c}_b^{\,2} \textrm{S}_B : \textrm{c}_c^{\,2} \textrm{S}_C], [ \textrm{c}_a^{\,2} \textrm{S}_A : - \textrm{S}_C \textrm{S}_A :\textrm{c}_c^{\,2} \textrm{S}_C],\\{[ \textrm{c}_a^{\,2} \textrm{S}_A :\textrm{c}_b^{\,2} \textrm{S}_B: - \textrm{S}_A \textrm{S}_B ]}$. 
\end{itemize}
Closely connected with it is the Simson cubic 
	\[
	\textrm{S}_A x_1(x_2^{2} + x_3^{2}) + \textrm{S}_B x_2(x_3^{2} + x_1^{2}) + \textrm{S}_C x_3(x_1^{2} + x_2^{2}) - 2 x_1 x_2 x_3 (1- \textrm{c}_a \textrm{c}_b \textrm{c}_c) = 0,
\]
which is the locus of tripoles of the Simson lines.\vspace*{-2mm}\\

The Euler-Feuerbach cubic belongs to a set of cubics which can be constructed as follows: Consider the pencil of circumconics of $\Delta$ which pass through a given point  $P = [p_1:p_2,p_3]$ different from $A, B, C$. The symmetry points of all these conics lie on a cubic which passes through the traces of $P$ and the traces of $G_i, i = 0,1,2,3.$ In metric affin geometries this cubic splits into the bicevian conic of $P$ and $G$ and the tripolar line of $G$.

\begin{figure}[!htbp]
\includegraphics[height=8cm]{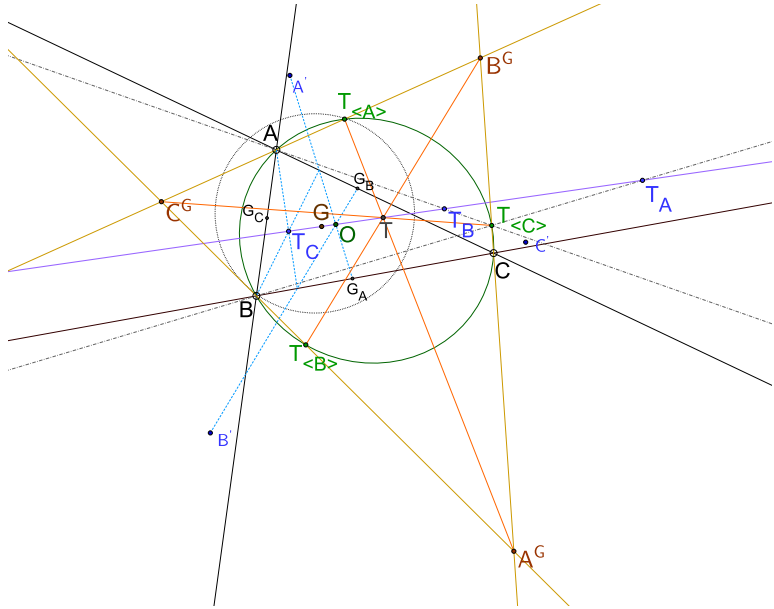}
\caption{The triplex points on the line $G \vee O$ and the points $T_{<A>}, T_{<B>},T_{<C>}$ on the circumcircle.}
\end{figure}

\subsubsection{} \textbf{Triplex points on $G \vee O$} \hspace*{\fill} \\
\noindent\hspace*{5mm}In euclidean geometry, triplex points were introduced by K. M\"utz [16]; further studies on triplex and related points have been carried out by E. Schmidt [20].
By joining the meet of the perpendicular bisector of $A\vee B$ and the side line $A\vee C$ with the vertex B and joining the meet of the perpendicular bisector of $A\vee C$ and the side line $A\vee B$ with the vertex C, we get two lines that meet at a point, the \textit{triplex point} $T_A$:  
	\[
\begin{split}			
	T_A &= (((C_G \vee C' ) \wedge (A \vee C))\vee B) \wedge (((B_G \vee B' ) \wedge (A \vee B))\vee C)\\
	&= [1: \frac {1-c_b}{c_a-c_c}:\frac {1-c_c}{c_a-c_b}].
\end{split}		
\]
The triplex points $T_B, T_C$ are defined accordingly.
It can be easily checked that $T_A, T_B, T_C$ are points on the line $G\vee O$. (See [16, 20] for the euclidean version.)\\
\hspace*{5mm}The points $T_{<A>}:= (B\vee T_C) \wedge (C\vee T_B), T_{<B>}:= (C\vee T_A) \wedge (A\vee T_C), T_{<C>}:= (A\vee T_B) \wedge (B\vee T_A)$ lie on the circumcircle; their coordinates are:
\begin{itemize}
\item[]$T_{<A>} = [1 - \textrm{c}_a: \textrm{c}_b - \textrm{c}_c: \textrm{c}_c - \textrm{c}_b],$
\item[]$T_{<B>} = [\textrm{c}_a - \textrm{c}_c: 1 - \textrm{c}_b: \textrm{c}_c - \textrm{c}_a],$
\item[]$T_{<C>} = [\textrm{c}_a - \textrm{c}_b: \textrm{c}_b - \textrm{c}_a: 1 - \textrm{c}_c].$
\end{itemize} 
Furthermore, these points lie on the lines $A\vee A^G, B\vee B^G, C\vee C^G$, respectively, and
the lines $A^G \vee T_{<A>}, B^G \vee T_{<B>},C^G \vee T_{<C>}$ meet at point
	\[ 
\begin{split}	
T :=\, &[- 3 \textrm{c}_a^2 + \textrm{c}_b^2 + \textrm{c}_c^2  - 2 \textrm{c}_b \textrm{c}_c + 2 \textrm{c}_a \textrm{c}_c + 2 \textrm{c}_a \textrm{c}_b  + 2 \textrm{c}_a  - 2 \textrm{c}_b- 2 \textrm{c}_c + 1:\cdots:\cdots]\\
=\, &[- 3 \textrm{s}_{a/2}^{\;4}+ \textrm{s}_{b/2}^{\;4} + \textrm{s}_{c/2}^{\;4} + 2 \textrm{s}_{a/2}^{\;2}\textrm{s}_{b/2}^{\;2} + 2 \textrm{s}_{a/2}^{\;2}\textrm{s}_{c/2}^{\;2}  - 2 \textrm{s}_{b/2}^{\;2} \textrm{s}_{c/2}^{\;2}:\cdots:\cdots]
\end{split}
\]
on the line $G \vee O$.\\\vspace*{-2mm}

\subsection{} \textbf{The line $O \vee K$} \hspace*{\fill} 

\begin{figure}[!htbp]
\includegraphics[height=8cm]{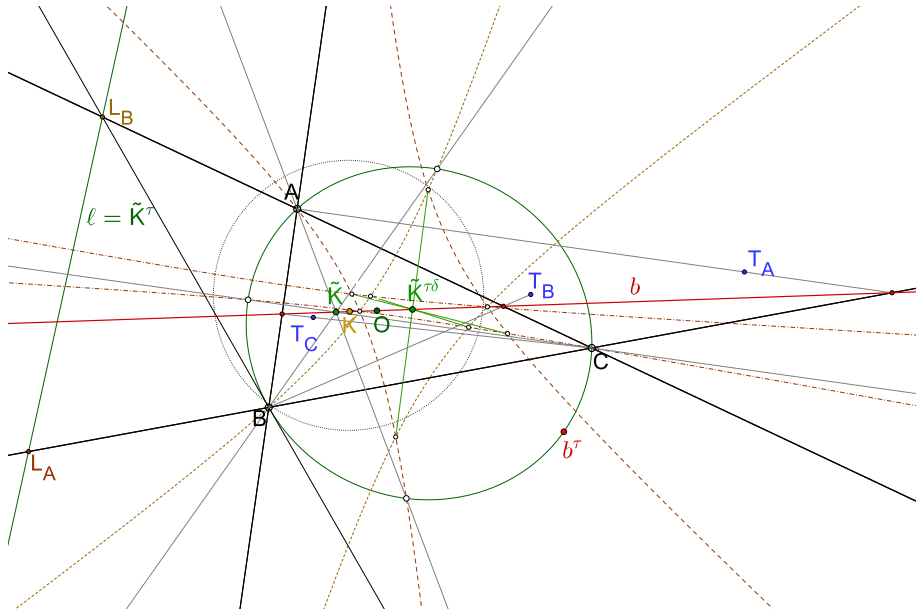}
\caption{The line $b = O \vee K$, the Lemoine axis $\ell = \tilde{K}^\tau$, the circumcircle and the apollonian circles.}
\end{figure}

\noindent\vspace*{-2mm}\subsubsection{} \textbf{Triangle centers on the line $O \vee K$} \hspace*{\fill} \\
\noindent\hspace*{5mm}The points $(A \vee T_A)\wedge (B\vee C), (B \vee T_B)\wedge (C\vee A), (C \vee T_C)\wedge (A\vee B)$ are collinear, they all lie on the line $O\vee K$.
The line $O \vee K$ is described by the equation 
	\[\frac{\textrm{c}_b - \textrm{c}_c}{1- \textrm{c}_a} x_1+ \frac{\textrm{c}_c - \textrm{c}_a}{1-\textrm{c}_b} x_2 + \frac{\textrm{c}_a - \textrm{c}_b}{1 - \textrm{c}_c}x_3=0.
\]
It has a tripole on the circumcircle and, besides the points already mentioned, it contains the Lemoine point $\tilde{K}$.\\
\hspace*{5mm}The \textit{Lemoine axis} $\tilde{K}^\tau$ is perpendicular to $O \vee K$, so its dual point also lies on $O \vee K$. This axis $\tilde{K}^\tau$ is also the polar of the Lemoine point with respect to the circumcircle and intersects the sidelines of $\Delta_0$ at the points
$L_A = [0:1{-}\textrm{c}_b:\textrm{c}_c{-}1],\; L_B = [\textrm{c}_a{-}1:0:1{-}\textrm{c}_c],\; L_C = [1{-}\textrm{c}_a:\textrm{c}_b{-}1:0].$\vspace*{-3 mm}\\
\subsubsection{} \textbf{The apollonian circles} \hspace*{\fill} \\
\noindent\hspace*{5mm}The line $O \vee K$ is the common radical axis of the circles $\mathcal{C}(L_A,A)$, $\mathcal{C}(L_B,B)$, $\mathcal{C}(L_C,C)$, which we will call \textit{apollonian circles} of $\Delta_0$.\\

\noindent\textit{Proof:} The two points 
$[(1-\textrm{c}_a)(1+\textrm{c}_a-t_{\pm}):(1-\textrm{c}_b)(1+\textrm{c}_b-t_{\pm}):(1-\textrm{c}_b)(1+\textrm{c}_b-t_{\pm})]$, \\ 
\hspace*{10mm} with $t_{\pm} =\frac{\sqrt{3}}{6} \big(\sqrt{1+ \textrm{c}_a + \textrm{c}_b + \textrm{c}_c} \pm \sqrt{1 - \textrm{c}_a^2 - \textrm{c}_b^2 - \textrm{c}_c^2 + 2 \textrm{c}_a \textrm{c}_b \textrm{c}_c}\big) $,\\
lie on $O \vee K$ and on each of the apollonian circles.$\;\;\Box$\vspace*{1 mm}\\ 
\noindent\textit{Remark:} The euclidean limits of these two points are the \textit{isodynamic points}. \vspace*{-3 mm}\\ 
\subsubsection{} \textbf{The Lemoine conic} \hspace*{\fill} \\
\noindent\hspace*{5mm}Define the point $P_1:= \textrm{par}(B \vee C,\tilde{K}) \wedge (B \vee C)$ and $P_2, P_3$ accordingly. Further, define $P_{23} := \textrm{par}(B \vee A,\tilde{K}) \wedge (A \vee C)$ and the points $P_{32}, P_{12},P_{21},P_{31},P_{13}$ accordingly. (Here we consider $A, B, C$ as the first, second and third point of the triangle $\Delta_0$,  respectively.)\\ The points $P_1, P_2, P_3$ lie on the line $\tilde{K}^\delta$ with the equation\vspace*{-1.5 mm}
	\[ (1{-}c_a{+}c_b{+}c_c{-}2 c_b c_c) x_1 + (1{+}c_a{-}c_b{+}c_c{-}2 c_c c_a) x_2{+}(1{+}c_a{+}c_b{-}c_c{-}2 c_a c_b) x_3 = 0. 
\]\vspace*{-2 mm}
The six points $P_{23},P_{32},P_{31},P_{13},P_{12}, P_{21} $ lie on a conic with the equation \vspace*{1 mm}
	\[ 
\begin{split}		
	\sum \limits_{cyclic}\!\Big(\big(&\nu_1{+}\nu_2{+}\nu_3{-}2 \nu_2\nu_3\big)\big(\nu_1(\nu_2{+}\nu_3{-}4\nu_2\nu_3){+}(\nu_2{-}\nu_3)^2\big)\nu_2 \nu_3 x_1^{\;2} \\
	- &\Big((\nu_1^4{+}\nu_1^3(3(\nu_2{+}\nu_3){-}8\nu_2\nu_3)){+}\nu_1^2(3(\nu_2^2{+}\nu_3^2){+}8\nu_2\nu_3{-}14\nu_2\nu_3(\nu_1{+}\nu_3){+}20\nu_2^2\nu_3^2)\\ 
	 \;\;\;\;\;& - \nu_1(\nu_2{+}\nu_3)(6\nu_2\nu_3(1{+}\nu_2{+}\nu_3){-}(\nu_2^2{+}\nu_3^2)){+}2\nu_2\nu_3(\nu_2{+}\nu_3)^2\Big)\nu_1 x_2 x_3\Big) \;\;=\;\;0,
\end{split}	
\]where we put $\nu_1 := 1-c_a = 2 s_{a/2}^{\;2}\,,\;\nu_2 := 1-c_b\,,\; \nu_3 := 1-c_c.$\\
We call this conic \textit{Lemoine conic}.\\
It can be proved by calculation:
\begin{itemize} 
\item[-] The line $K \vee O$ is a symmetry line of the Lemoine conic. 
\item[-] If the line  $\tilde{K}^\delta$ has common points with the circumcircle, then these points are also points on the Lemoine conic.
\item[-] The pole of $\tilde{K}^\delta$ with respect to the circumcircle is a point on $K \vee O$.
\item[-] The pole of $\tilde{K}^\delta$ with respect to the Lemoine conic is a point on $K \vee O$.
\end{itemize}  \vspace*{15 mm}

\subsection{} \textbf{The Akopyan line $O \vee H^{\star}$}\hspace*{\fill} 

\begin{figure}[!htbp]
\includegraphics[height=8cm]{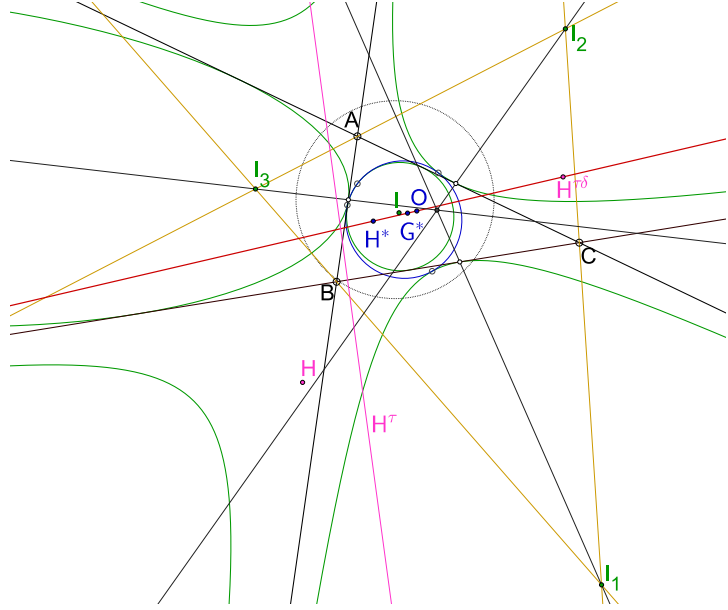}
\caption{The Akopyan line and the Hart circle together with the incircle and the excircles of $\Delta_0$.}
\end{figure}
\subsubsection{} \textbf{Triangle centers on the line $O \vee H^{\star}$}\hspace*{\fill} \\
\noindent\hspace*{5mm}There are several triangles centers, introduced by Akopyan [2], lying on the join of the circumcenter $O$ and the orthostar $H^\star$; therefore, we propose to name this line \textit{Akopyan line}.
(Akopyan uses the name \textit{Euler line} for it, Vigara the name \textit{Akopyan Euler line}.) Its equation is
	\[ 
	\begin{split}	
	(c_b{-}c_c) (1{+}2 \textrm{c}_a{-}\textrm{c}_b{-}\textrm{c}_c{-}\textrm{c}_b \textrm{c}_c) x_1 &+ (\textrm{c}_c{-}\textrm{c}_a) (1{-}\textrm{c}_a{+}2 \textrm{c}_b{-}\textrm{c}_c{-}\textrm{c}_c \textrm{c}_a) x_2 \\
	&+ (\textrm{c}_a{-}\textrm{c}_b) (1{-}\textrm{c}_a{-}\textrm{c}_b{+}2 \textrm{c}_c{-}\textrm{c}_a \textrm{c}_b)x_3 = 0.
\end{split}	
\]
 \\
\noindent\hspace*{5mm}As a first point on this line, apart from $O$ and $H^{\star}$, we introduce the point $G^\sharp$, whose cevians bisect the triangle area in equal parts. The existence of such a point was already shown by J. Steiner [21]. (See also [3].) The coordinates of $G^\sharp$ are
	\[ 
\begin{split}		
	&[\frac {\sqrt{1{+}\textrm{c}_a}}{\sqrt{2}\,\sqrt{1{+}\textrm{c}_a}+\sqrt{1{+}\textrm{c}_b}\,\sqrt{1{+}\textrm{c}_c}}:\frac {\sqrt{1{+}\textrm{c}_b}}{\sqrt{2}\,\sqrt{1{+}\textrm{c}_b}+\sqrt{1{+}\textrm{c}_c}\,\sqrt{1{+}\textrm{c}_a}}\\
	&\hspace*{45mm}:\frac {\sqrt{1{+}\textrm{c}_c}}{\sqrt{2}\,\sqrt{1{+}\textrm{c}_c}+ \sqrt{1{+}\textrm{c}_a}\,\sqrt{1{+}\textrm{c}_b}}] \\
	= \;&[\frac {\textrm{c}_{a/2}}{\textrm{c}_{a/2}+\textrm{c}_{b/2}\textrm{c}_{c/2}}:\frac {\textrm{c}_{b/2}}{\textrm{c}_{b/2}+\textrm{c}_{c/2}\textrm{c}_{a/2}}:
	\frac {\textrm{c}_{c/2}}{\textrm{c}_{c/2}+\textrm{c}_{a/2}\textrm{c}_{b/2}}].
\end{split}	
\]
The calculation is carried out according to the construction of $G^\sharp$ described below.\\

Akopyan [2] shows that the cyclocevian of $G^\sharp$ lies on the line $O \vee G^\sharp$ and has properties that justify to call it a \textit{pseudo-orthocenter}. We denote this point by $H^\sharp$; its coordinates are\vspace*{-1mm}\\
\[H^\sharp = [\frac {\textrm{c}_{a/2}}{\textrm{c}_{a/2}-\textrm{c}_{b/2}\textrm{c}_{c/2}}:\frac {\textrm{c}_{b/2}}{\textrm{c}_{b/2}-\textrm{c}_{c/2}\textrm{c}_{a/2}}:\frac {\textrm{c}_{c/2}}{\textrm{c}_{c/2}-\textrm{c}_{a/2}\textrm{c}_{b/2}}].
\]
The vertices $A$ and $B$ together with the traces $A_{H^\sharp}, B_{H^\sharp}$ of $H^\sharp$ lie on a circle with center \vspace*{-4mm}\\
\[
[\frac {\textrm{c}_{a/2}}{\textrm{c}_{a/2}+\textrm{c}_{b/2}\textrm{c}_{c/2}}:\frac {\textrm{c}_{b/2}}{\textrm{c}_{b/2}+\textrm{c}_{c/2}\textrm{c}_{a/2}}: \frac {\textrm{c}_{c/2} \textrm{s}_{c/2}^2}{(\textrm{c}_{a/2} + \textrm{c}_{b/2} \textrm{c}_{c/2}) (\textrm{c}_{b/2}+\textrm{c}_{a/2} \textrm{c}_{c/2})}].
\]
This center is also the intersection of the perpendicular bisector of $[A,B]_+$ and the cevian of $G^\sharp$ through the vertex $C$.\\
(Mutatis mutandis, we can replace $A, B, A_{H^\sharp}, B_{H^\sharp}$ by $B, C, B_{H^\sharp}, C_{H^\sharp}$ or by $C, A,$ $C_{H^\sharp}, A_{H^\sharp}$.)\vspace*{1mm}\\
Akopyan also shows that the center $N^\sharp$ of the common cevian circle of $G^\sharp$ and $H^\sharp$ is a point on $O \vee G^\sharp$. The barycentric re\-pre\-sen\-tation of this point is \vspace*{-1mm}
	\[ 
\begin{split}
N^\sharp &= [(\textrm{c}_a{+}1)(\textrm{c}_a(\textrm{c}_b\,\textrm{c}_c{-}1) + 1 - \textrm{c}_b^2 + \textrm{c}_b\,\textrm{c}_c - \textrm{c}_c^2):\cdots:\cdots]\\
	&= [\textrm{c}_{a/2}^{\;2} (\textrm{c}_{a/2}^{\;2} (2\textrm{c}_{b/2}^{\;2}\textrm{c}_{c/2}^{\;2} -\textrm{c}_{b/2}^{\;2}-\textrm{c}_{c/2}^{\;2}) + \textrm{c}_{b/2}^{\;2} \textrm{s}_{b/2}^{\;2} + \textrm{c}_{c/2}^{\;2} \textrm{s}_{c/2}^{\;2}:\cdots:\cdots].\vspace*{2mm}		
\end{split}	
\] 
The points $H^\sharp,N^\sharp,G^\sharp,O$ form a harmonic range.\vspace*{1mm}\\
\textit{Proof\vspace*{-2mm}}: 
\[
\begin{split}		
\textrm{Put}\; \boldsymbol{h} := &(\frac {\textrm{c}_{a/2}}{\textrm{c}_{a/2}-\textrm{c}_{b/2}\textrm{c}_{c/2}},\frac {\textrm{c}_{b/2}}{\textrm{c}_{b/2}-\textrm{c}_{c/2}\textrm{c}_{a/2}},\frac {\textrm{c}_{c/2}}{\textrm{c}_{c/2}-\textrm{c}_{a/2}\textrm{c}_{b/2}}),\\
\boldsymbol{n} := &(\textrm{c}_{a/2}^{\;2} (\textrm{c}_{a/2}^{\;2} (2\textrm{c}_{b/2}^{\;2}\textrm{c}_{c/2}^{\;2} -\textrm{c}_{b/2}^{\;2}-\textrm{c}_{c/2}^{\;2}) + \textrm{c}_{b/2}^{\;2} \textrm{s}_{b/2}^{\;2} + \textrm{c}_{c/2}^{\;2} \textrm{s}_{c/2}^{\;2}), \\
                &\;\,\textrm{c}_{b/2}^{\;2} (\textrm{c}_{b/2}^{\;2} (2\textrm{c}_{c/2}^{\;2}\textrm{c}_{a/2}^{\;2} -\textrm{c}_{c/2}^{\;2}-\textrm{c}_{a/2}^{\;2}) + \textrm{c}_{c/2}^{\;2} \textrm{s}_{c/2}^{\;2} + \textrm{c}_{b/2}^{\;2} \textrm{s}_{b/2}^{\;2}),\\
								&\;\,\textrm{c}_{c/2}^{\;2} (\textrm{c}_{c/2}^{\;2} (2\textrm{c}_{a/2}^{\;2}\textrm{c}_{b/2}^{\;2} -\textrm{c}_{a/2}^{\;2}-\textrm{c}_{b/2}^{\;2}) + \textrm{c}_{a/2}^{\;2} \textrm{s}_{a/2}^{\;2} + \textrm{c}_{a/2}^{\;2} \textrm{s}_{a/2}^{\;2})),\\
\boldsymbol{g} := &(\frac {\textrm{c}_{a/2}}{\textrm{c}_{a/2}+\textrm{c}_{b/2}\textrm{c}_{c/2}},\frac {\textrm{c}_{b/2}}{\textrm{c}_{b/2}+\textrm{c}_{c/2}\textrm{c}_{a/2}},\frac {\textrm{c}_{c/2}}{\textrm{c}_{c/2}+\textrm{c}_{a/2}\textrm{c}_{b/2}}),\\
	\boldsymbol{o} := &(\textrm{s}_{a/2}^{\;2}(1{+}\textrm{c}_{a/2}^{\;2}{-}\textrm{c}_{b/2}^{\;2}{-}\textrm{c}_{c/2}^{\;2}) , \textrm{s}_{b/2}^{\;2}(1{-}\textrm{c}_{a/2}^{\;2}{+}\textrm{c}_{b/2}^{\;2}{-}\textrm{c}_{c/2}^{\;2}), \textrm{s}_{c/2}^{\;2}(1{-}\textrm{c}_{a/2}^{\;2}{-}\textrm{c}_{b/2}^{\;2}{+}\textrm{c}_{c/2}^{\;2})),\\
	r:=\, &\textrm{c}_{a/2} \textrm{c}_{b/2} \textrm{c}_{c/2},\\
	s:=\, &(\textrm{c}_{a/2}{+}\textrm{c}_{b/2} \textrm{c}_{c/2})(\textrm{c}_{b/2}{+}\textrm{c}_{c/2} \textrm{c}_{a/2})(\textrm{c}_{c/2}{+}\textrm{c}_{a/2} \textrm{c}_{b/2})(2r{+}1{-}\textrm{c}_{a/2}^{\;2}{-}\textrm{c}_{b/2}^{\;2}{-}\textrm{c}_{c/2}^{\;2}),\\
	t:=\, &(\textrm{c}_{a/2}{-}\textrm{c}_{b/2} \textrm{c}_{c/2})(\textrm{c}_{b/2}{-}\textrm{c}_{c/2} \textrm{c}_{a/2})(\textrm{c}_{c/2}{-}\textrm{c}_{a/2} \textrm{c}_{b/2})(2r{-}1{+}\textrm{c}_{a/2}^{\;2}{+}\textrm{c}_{b/2}^{\;2}{+}\textrm{c}_{c/2}^{\;2}),\vspace*{-2mm}
	\end{split}	
\]
and verify the equations \vspace*{-2mm}
\[
r \boldsymbol{o} + \boldsymbol{n} = s \boldsymbol{g} \;\;\textrm{and} \;\,r \boldsymbol{o} - \boldsymbol{n} = t \boldsymbol{h}.\vspace*{-1mm}
\]
The cevian circle of $G^\sharp$ can be seen as a good substitute in elliptic geometry for the euclidean nine point circle, even more so, since this circle, as also proved by Akopyan, touches the incircles of $\Delta_i$ for $i = 0,1,2,3$. The common cevian circle of $G^\sharp$ and $H^\sharp$ we like to name \textit{Hart circle} of $\Delta_0$, because A.\,S. Hart [9] calculated 1861 the equation of the circle which touches incircle and the excircles of a spherical triangle, and the name \textit{Hart's circle} is used by G. Salmon in [19]. Salmon showed that its center $N^\sharp$ lies on the lines $G \vee H$ and $O^+ \vee H^-$ and he calculated the (trilinear) coordinates of $N^\sharp$.\\
The Feuerbach point (touchpoint of the incircle with the Hart circle) is \vspace*{-1.5 mm}
\[F\!e = [\textrm{S}^2 + \textrm{S}_B \textrm{S}_C - \textrm{s}_a^{ 2} \textrm{s}_b \textrm{s}_c:\textrm{S}^2 + \textrm{S}_C \textrm{S}_A - \textrm{s}_a \textrm{s}_b^{ 2} \textrm{s}_c:\textrm{S}^2 + \textrm{S}_A \textrm{S}_B - \textrm{s}_a \textrm{s}_b \textrm{s}_c^{ 2}].\vspace*{1 mm}
\]
\noindent\textit{Remark}: The Akopyan line is a line which contains a point together with its circumcevian conjugate and the center of their common cevian circle. Such a line is called \textit{cevian axis} [6]. The orthoaxis and the line $G \vee O$ are, in general, not cevian axes of the triangle $\Delta_0$. But the line $G \vee O$ is a cevian axis of the anticevian triangle $(A^GB^GC^G)_0$ of $G$, see F\textsc{igure} 2. \\
\begin{figure}[!htbp]
\includegraphics[height=7cm]{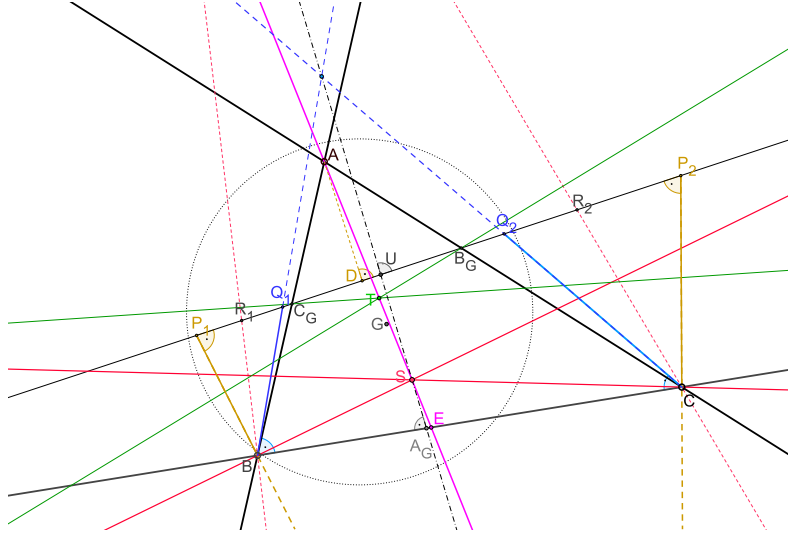}
\caption{Construction of the bisector $A \vee S$ through the vertex $A$.} 
\end{figure}

Explanations to F\textsc{igure} 5:  By projecting the points $A, B, C$ onto the sideline $B_G \vee C_G$ of the medial triangle of $\Delta_0$, we get the points $D, P_1$ and $P_2$, respectively. The area $2\epsilon$ of the triangle $\Delta_0$ is the same as the area of the quadrangle $(AP_1P_2B)_0$ because triangle $(C_GAD)_0$ is congruent to triangle $(C_GBP_1)_0$ and triangle $(B_GAD)_0$ congruent to triangle $(B_GBP_2)_0$. \\
The lines $\textrm{perp}(B \vee C,B)$ and  $\textrm{perp}(B \vee C,C)$ meet $B_G \vee C_G$ at $Q_1$ and $Q_2$, respectively. It follows that $2\epsilon = \mu(\angle_{+} P_1BQ_1) + \mu(\angle_{+} P_2CQ_2)$.\\ 
Let $W$ be the meet of the lines $B \vee C$ and $B_G \vee C_G$.  Confirm by calculation that the mirror image of $P_1BQ_1$ in $W$ (or in its dual line $W^\delta = A_G \vee U$) is $P_2BQ_2$. It follows that $\mu(\angle_{+} P_1BQ_1) = \mu(\angle_{+} P_2CQ_2)$ = $\epsilon$.\\
The lines $B \vee R_1$ and $C \vee R_2$ are the internal bisectors, the lines $B \vee S$ and $C \vee S$ the external bisectors of $\angle_{+} P_1BQ_1$ resp. $\angle_{+} P_2CQ_2$. Define E := $(A \vee S) \wedge (B \vee C)$ and let $T$ be the midpoint of $[A,E]_+$. We can now confirm by calculating that $T \vee C_G = perp (B \vee R_1, C_G)$ and conclude that the area of the triangle $(ABE)_0$ is $\epsilon$.\\

List of triangle centers:\vspace*{2mm} 

\noindent
\begin{tabular}{|c|c|c|}
   \hline   & & \\
   \!triangle center\!& $P=$ &\!euclidean limit point\!\\ 
   $P$ &$\!\![f(\alpha,\beta,\gamma){:}f(\beta,\gamma,\alpha){:}f(\gamma,\alpha,\beta)]$,\!\! & (We adopt the nota- \\ 
   & $f(\alpha,\beta,\gamma) =$ & tion from [11].)\\ 
	 & & \\ \hline \hline
   $\!\,\;\;\;\;\;\;\;\;\;\;G\;\;\;\;\;\;\;\;\;\;\,$ &\;\;\;\;\;\;\;\;\;\;\;\;\;\;\;\;\;\;\;\;\;\;\;\;\;\;\;\;\;\;\;\;\;\;\,1\;\;\;\;\;\;\;\;\;\;\;\;\;\;\;\;\;\;\;\;\;\;\;\;\;\;\;\;\;\;\;\;\;\;\,&  \\ \cline{1-2}
	 \rule{0pt}{10pt}\,\;\;$G^+\;\;$  &$1- \frac{2\,\sin \alpha \, \sin \alpha\text{-}\epsilon}{\sin \beta  \sin \gamma}$&  $X_2$\\ \cline{1-2}
	 \rule{0pt}{10pt}$\T G^{\sharp^{}}$ &$ \frac{\sin \alpha} {\sin \alpha\, + \,\sin {\alpha\text{-}{\frac \epsilon 2_{}}}}$ &\\ \hline
  $\!I$ & $\sin \alpha$ & $X_1$ \\\hline
	 $\!O\,$ & $\sin \alpha \;\,\cos \alpha{\text{-}}\epsilon$  &\multirow{2}*{$X_3$}\\ \cline{1-2}
	 \rule{0pt}{9pt}\,\T $O^+$ & $\sin 2\alpha$ &  \\\hline
 $\!H\,$\, & $\tan \alpha$   &\T \multirow{3}*{$X_4$}  \\ \cline{1-2}
	 $\T\;H^-$\, &\!$\sin \alpha \,/ \cos \alpha{\text{-}}\epsilon$\! & \\ \cline{1-2}
	$\T H^\sharp$\, &$ \frac{\sin \alpha} {\sin \alpha\,  - \,\sin {\alpha\text{-}{\frac \epsilon 2_{}}}}$& \\\hline
	$\T \!N^\sharp$ & $\sin \alpha\;\,\cos \beta{\text{-}}\gamma$ & $X_5$ \\\hline
	\rule{0pt}{3pt}$\!K$\, & $\sin^{\scriptscriptstyle{2}} \alpha$ &\multirow{2}*{$X_6$}  \\ \cline{1-2}
	\rule{0pt}{11pt}$\!\tilde{K}$\, & $\sin \alpha  \,\sin \alpha{\text{-}}\epsilon$&  \\ \hline
		 $Ge$\, & $\tan {\scriptstyle{\frac 1 2}} \alpha$ & $X_7$ \\ \hline
	 $N\!a$\, & $\cot {\scriptstyle{\frac 1 2}} \alpha$ & $X_8$ \\ \hline
	$\T F\!e$\, & $\sin \alpha - \sin\alpha\, \cos \beta{\text{-}}\gamma$ & $X_{11}$ \\ \hline
	\rule{0pt}{10pt} \multirow{2}*{$\!L\;\;\,$}  &$\!3\xi_{\alpha}^{\;2}{-}2\xi_{\alpha}(\xi_{\beta}{+}\xi_{\gamma}){-}(\xi_{\beta}{-}\xi_{\gamma})^2{-}\phi_\alpha(\xi_\alpha^{\;2}{-}\xi_\beta^{\;2}{-}\xi_\gamma^{\;2}),\!$  &\multirow{3}*{$X_{20}$} \\ 
		 \rule{0pt}{10pt}&with $\xi_\alpha = \textrm{s}_{\alpha} \textrm{s}_{\alpha{\text{-}\epsilon}}, \dots$ and $\phi_\alpha =\frac{2 \sin \epsilon \, \sin \alpha\text{-}\epsilon}{\sin \beta  \sin \gamma}$   &   \\  \cline{1-2}
	 \rule{0pt}{10pt} $\;\;T\;$(cf. 3.2.3.) &$-3\textrm{s}_{\alpha/2}^{\;4} + 2\textrm{s}_{\alpha/2}^{\;2}(\textrm{s}_{\beta/2}^{\;2}{+}\textrm{s}_{\gamma/2}^{\;2}) + (\textrm{s}_{\beta/2}^{\;2}{-}\textrm{s}_{\gamma/2}^{\;2})^2  $&\\ \hline
\rule{0pt}{10pt}$\;\;\;\;\;\;\;\;H^{\tau \delta}\;\;\;\;\;\;\;$  &$\sin \alpha \,(\cos \alpha - 2 \cos \beta \;\cos \gamma)$&{$X_{30}$}\\ \hline
\multirow{2}*{$(O \vee K)^{\tau}$ } &$\xi_\alpha/(\xi_{\beta}{-}\xi_{\gamma}) $  &\multirow{2}*{$X_{110}$}  \\ 
\rule{0pt}{10pt}& $\;\;\;\;\;\;\xi_\alpha = \textrm{s}_{\alpha} \textrm{s}_{\epsilon{\text{-}}\alpha}$, $\xi_\beta = \textrm{s}_{\beta} \textrm{s}_{\epsilon{\text{-}}\beta}$, $\xi_\gamma = \textrm{s}_{\gamma} \textrm{s}_{\epsilon{\text{-}}\gamma}\;\;\;\;\;\;$   & \\  \hline
\rule{0pt}{10pt}\multirow{2}*{\,\;$\tilde{K}^{\tau \delta}$ } &$\xi_\alpha(\xi_{\beta}{+}\xi_{\gamma}){-}\xi_{\beta}^{\;2}{-}\xi_{\gamma}^{\;2} - 2\,\textrm{s}_\epsilon\, \textrm{s}_{\alpha\text{-}\epsilon}\, \textrm{s}_{\beta\text{-}\epsilon}\, \textrm{s}_{\gamma\text{-}\epsilon}, $  &{infinity point }  \\ 
\rule{0pt}{10pt}& $\xi_\alpha = \textrm{s}_{\alpha} \textrm{s}_{\alpha\text{\textrm-}\epsilon}$, $\xi_\beta = \textrm{s}_{\beta} \textrm{s}_{\beta{\text{-}}\epsilon}$, $\xi_\gamma = \textrm{s}_{\gamma} \textrm{s}_{\gamma{\text{-}\epsilon}}$   &\,on the Brocard axis\, \\  \hline
\end{tabular}\\

\end{document}